\documentstyle[12pt]{article}  
\def\sq{\hbox {\rlap{$\sqcap$}$\sqcup$}}
\def\sq{\hbox {\rlap{$\sqcap$}$\sqcup$}}
\def\R{ {\rm R \kern -.31cm I \kern .15cm}}
\def\C{ {\rm C \kern -.15cm \vrule width.5pt \kern .12cm}}
\def\Z{ {\rm Z \kern -.27cm \angle \kern .02cm}}
\def\N{ {\rm N \kern -.26cm \vrule width.4pt \kern .10cm}}
\def\1{{\rm 1\mskip-4.5mu l} }
\def\lsim{\raise0.3ex\hbox{$<$\kern-0.75em\raise-1.1ex\hbox{$\sim$}}}
\def\gsim{\raise0.3ex\hbox{$>$\kern-0.75em\raise-1.1ex\hbox{$\sim$}}}
\def\noi{\noindent}

\def\beq{\begin{equation}}   \def\eeq{\end{equation}}
\def\bea{\begin{eqnarray}}  \def\eea{\end{eqnarray}}
\def\nn{\nonumber}
\def\noi{\noindent}
\def\beeq{\begin{eqnarray}} \def\eeeq{\end{eqnarray}}
\newcommand\mysection{\setcounter{equation}{0}\section}

\newcounter{hran}

\begin{document} 
\centerline{\large\bf Long range scattering for some }
 \vskip 3 truemm \centerline{\large\bf  Schr\"odinger related nonlinear systems\footnote{Expanded version of a lecture given by J. G. at the COE Symposium on Nonlinear Dispersive Equations, Sapporo, September 2004.}}  

\vskip 0.5 truecm

\centerline{\bf J. Ginibre}
\centerline{Laboratoire de Physique Th\'eorique\footnote{Unit\'e Mixte de
Recherche (CNRS) UMR 8627}}  \centerline{Universit\'e de Paris XI, B\^atiment
210, F-91405 ORSAY Cedex, France}
\vskip 3 truemm
\centerline{\bf G. Velo}
\centerline{Dipartimento di Fisica, Universit\`a di Bologna}  \centerline{and INFN, Sezione di
Bologna, Italy}

\vskip 1 truecm

\begin{abstract} We review some results, especially recent ones, on the
theory of scatttering and more precisely on the local Cauchy problem at
infinity in time, in long range situations, for some nonlinear systems
including some form of the Schr\"odinger equation. We consider in
particular the Wave-Schr\"odinger system (WS)$_3$ and the
Maxwell-Schr\"odinger system (MS)$_3$ in space dimension $n=3$ and the
Klein-Gordon-Schr\"odinger system (KGS)$_2$ in space dimension $n=2$.
We also consider the Zakharov system (Z)$_n$ which can be studied from
the same point of view in space dimension $n=3$ where it is short range
and in space dimension $n=2$ where it is in a mixed situation. We
concentrate on the applications of a direct method which is
intrinsically restricted to the case of small Schr\"odinger data and to
the borderline long range case, to which (WS)$_3$, (MS)$_3$ and
(KGS)$_2$ belong. The main results in all cases are the existence of
solutions defined for large times and with prescribed asymptotic
behaviour, without any size restriction for the data of the wave,
Maxwell or Klein-Gordon function. Furthermore convergence rates as
negative powers of $t$ are obtained for the solutions in suitable norms.
The asymptotic forms in the long range cases are obtained by suitable
modifications of the solutions of the underlying free linear system.
\end{abstract}

\vskip 1 truecm
\noi AMS Classification : Primary 35P25. Secondary 35B40, 35Q40, 81U99.  \par \vskip 2 truemm

\noi Key words : Long range scattering, Wave-Schr\"odinger system, Maxwell-Schr\"odinger system, Zakharov system, Klein-Gordon-Schr\"odinger system.\par 
\vskip 1 truecm

\noindent LPT Orsay 04-122\par
\noindent November 2004\par

\newpage
\pagestyle{plain}
\baselineskip 18pt
\mysection{Introduction}
\hspace*{\parindent} The purpose of this paper is to review some results, especially recent
ones, on the theory of scattering in long range situations for some non
linear systems containing some form of the Schr\"odinger equation. In
order to put the subject in perspective, we first list a number of
equations and systems for which one can ask the same questions and to
which one can apply the same methods as described below. In all this
paper, $u$ denotes a complex valued function defined in space time
${I\hskip-1truemm R}^{n+1}$, $\Delta$ denotes the Laplacian in ${I\hskip-1truemm R}^n$ and $\sq = \partial_t^2
- \Delta$ the d'Alembertian in ${I\hskip-1truemm R}^{n+1}$. When appended to the symbol of
an equation or system, the subscript $n$ means that that equation or
system is considered in space dimension $n$. \par

The starting point is the linear Schr\"odinger equation 
$$i \partial_t u = - (1/2) \Delta u + V u \eqno({\rm LS)_n}$$

\noi where $V$ is a real potential defined in ${I\hskip-1truemm R}^n$,
a typical form of which is
\beq
\label{1.1e}
V(x) = \kappa |x|^{-\gamma}
\eeq

\noi with $\kappa \in {I\hskip-1truemm R}$ and $\gamma > 0$. Next come
the nonlinear Schr\"odinger equation, a simple form of which is
$$i \partial_t u = - (1/2) \Delta u + \kappa |u|^{p-1} u \eqno({\rm NLS)_n}$$

\noi with $\kappa \in {I\hskip-1truemm R}$, with $1 < p < \infty$ for $n=1,2$ and $1 < p < (n+2)/(n-2)$ for $n \geq 3$, and the Hartree equation
$$i \partial_t u = - (1/2) \Delta u + (V \ *\  |u|^2)u \eqno({\rm R3)_n}$$

\noi where $V$ is a real even potential defined in ${I\hskip-1truemm
R}^n$ of the same type as for (LS)$_n$, for instance given by
(\ref{1.1e}), and where $*$ denotes the convolution in
${I\hskip-1truemm R}^n$. Next come some systems including the
Schr\"odinger equation with a real time dependent potential $A$
defined in ${I\hskip-1truemm R}^{n+1}$ and a second equation whereby
the potential $A$ is nonlinearly coupled to the Schr\"odinger function.
These include the Wave-Schr\"odinger system 
$$\left \{ \begin{array}{l} i \partial_t u = - (1/2) \Delta u + A u\\ \\\sq A = - |u|^2\end{array} \right . \eqno({\rm WS)_n}$$

\noi which will be of special interest for $n=3$, the Klein-Gordon-Schr\"odinger system
$$\left \{ \begin{array}{l} i \partial_t u = - (1/2) \Delta u + A u\\ \\(\sq + 1) A = - |u|^2\end{array} \right . \eqno({\rm KGS})_{\rm n}$$

\noi which will be of special interest for $n=2$, and the Zakharov system, the simplest version of which is
$$\left \{ \begin{array}{l} i \partial_t u = - (1/2) \Delta u + A u\\ \\\sq  A = \Delta |u|^2\end{array} \right . \eqno({\rm Z})_{\rm n}$$

\noi which will be of special interest for $n = 2,3$. Finally, we shall
consider the Maxwell-Schr\"odinger system in space dimension 3, which
we write here in the Coulomb gauge
$$\left \{ \begin{array}{l} i \partial_t u = - (1/2) \Delta_A u + g(|u|^2) u\\ \\\sq A = P\ {\rm Im}\ \overline{u} \nabla_A u\ , \quad \nabla \cdot A = 0\ . \end{array} \right . \eqno({\rm MS)_3}$$

\noi Here $A$ is an ${I\hskip-1truemm R}^3$ valued function (the
magnetic potential) defined in ${I\hskip-1truemm R}^{3+1}$, $\nabla_A = \nabla - iA$
and $\Delta_A = \nabla_A^2$ are respectively the covariant gradient and
the covariant Laplacian, $g(|u|^2)$ is the Hartree
interaction
\beq
\label{1.2e}
g(|u|^2) = (4 \pi |x|)^{-1} \ * \ |u|^2
\eeq 

\noi which is of the same type as in (R3)$_3$ with $V$ given by
(\ref{1.1e}) with $\kappa = (4 \pi)^{-1}$ and $\gamma = 1$, and $P$ is
the projector on divergence free vector fields
\beq
\label{1.3e}
P = \1 - \nabla \Delta^{-1} \nabla \ .
\eeq

\noi The condition $\nabla \cdot A = 0$ is the Coulomb gauge condition.
We shall describe the (MS)$_3$ system in more detail in Section~3.\par

All the previous equations and systems are Lagrangian, namely are the
Euler-Lagrange equations associated with some suitable Lagrange
function. As a consequence of that fact and of some obvious invariance
properties, they possess some formally conserved quantities. In
particular the $L^2$ norm of $u$ is formally conserved in all cases, as
well as an energy function which is closely related to and in favorable
cases controls the $H^1$ norm of $u$ and some norm of $A$. This leads
to the definition of an energy space which is $H^1$ for $u$ for
(LS)$_n$, (NLS)$_n$ and (R3)$_n$, which is $H^1 \oplus H^1 \oplus L^2$
for $(u, A, \partial_t A)$ for (WS)$_n$, (KGS)$_n$ and (MS)$_3$ and
which is $H^1 \oplus L^2\oplus \dot{H}^{-1}$ for $(u, A, \partial_t A)$
for (Z)$_n$. We shall occasionally refer to those energy spaces in the
sequel, but we shall never make use of the explicit form of the
energy.\par

We now review briefly some basic ideas of scattering theory in order to
extract therefrom the mathematical problem to be addressed in this
paper. In a typical scattering experiment, the experimentalist has at
his/her disposal a target at rest, typically in some stationary state,
and does his/her best to prepare an incoming beam of free waves or free
particles. Thus the mathematical picture is that of a dynamical system
target $\times$ beam with an asymptotic motion in the distant past of a
simple nature, typically stationary state $\times$ free motion. As the
experiments proceeds, the beam reaches the target, thereby giving rise
to a complicated interacting motion. Eventually a modified more or less
free beam emerges, possibly leaving the target in a different state, so
that the asymptotic motion in a distant future is again simple. The
basic objects describing that process are the wave operators $\Omega_{\pm}$
for the past $(-)$ and future $(+)$, defined as the maps going from the
asymptotic simple motion in the past and/or future to the actual
interacting motion. The experimentalist then measures the scattering
operator $S = \Omega_+^{-1} \circ \Omega_-$, namely the map from the
simple asymptotic motion in the past to the simple asymptotic motion in
the future. The basic idea that we shall keep from this picture is that
we have a dynamical system with simple asymptotics in time, and that we
want to classify the possible motions of the system through their
asymptotics by constructing the wave operators.\par

We now give a slightly different and more formal description of the
previous picture, taking the simple case of one single
Schr\"odinger-like equation, e.g. (LS)$_n$, (NLS)$_n$ or (R3)$_n$. The
dynamical system has the space of Cauchy data as state space, and the
evolution is described by the equation at hand. The extension to
systems of two equations is obvious. We give ourselves a set ${\cal
U}_a$ of presumed asymptotic motions $u_a$ for the system at hand,
parametrized by some data $u_+$. In the present case, a natural
candidate for ${\cal U}_a$ would be the set of solutions of the free
Schr\"odinger equation, namely
\beq
\label{1.4e}
u_a (t) = U(t) u_+
\eeq

\noi where $U(t)$ is the unitary group
\beq
\label{1.5e}
U(t) = \exp (i(t/2)\Delta )
\eeq

\noi solving the free Schr\"odinger equation. The initial data $u_+$ is
called the asymptotic state. One is then faced with two natural
problems. \\

\noi \underbar{Problem 1.} Given $u_a \in {\cal U}_a$, construct a
solution $u$ of the equation at hand such that $u(t) - u_a(t)$ tends to
zero as $t \to + \infty$ in a suitable sense, more precisely in
suitable norms. If that problem can be solved for ``all'' $u_a$, one can
define the wave operator $\Omega_+$ (for positive times) as the map
$u_a \to u$ thereby obtained. The same problem can be considered for $t
\to - \infty$, thereby giving rise to the wave operator $\Omega_-$ for
negative times. Thus Problem 1 is that of the construction of the wave
operators.\par

Traditionally the wave operator, for instance $\Omega_+$ is defined as
the map $u_+ \to u(0)$ where $u$ is the solution considered above. Thus
the problem decomposes into two steps. The first step is to construct
the required solution $u$ in a neighborhood of infinity in time,
namely in an interval $[T, \infty )$ for $T$ sufficiently large and
with given asymptotic behaviour $u_a$. This is the local Cauchy problem
at infinity in time. The second and rather independent step consists in
extending the solution previously obtained down from $t = T$ to $t=0$,
and therefore reduces to the global Cauchy problem at finite times. In
this paper we shall concentrate on the first step and leave aside the
second one. Actually the latter is well controlled, for instance in the
energy space, for all systems in the previous list for the relevant
values of $n$, except for the (MS)$_3$ system (see below in Section 3
for that system).\\

\noi \underbar{Problem 2.} This is the converse to Problem 1. Given a
generic solution $u$ of the equation at hand, find an asymptotic motion
$u_a \in {\cal U}_a$ such that $u(t) - u_a(t)$ tends to zero as $t \to
+ \infty$ in a suitable sense (in suitable norms). If that problem and
the same one for $t \to - \infty$ can be solved (in a suitable
functional framework) for ``all '' $u$, one says that asymptotic
completeness holds with respect to the set ${\cal U}_a$. Asymptotic 
completeness is a much harder question than the existence of the wave
operators. It requires in particular that all possible asymptotic
behaviours of the solutions of the equation at hand have been
identified and included in ${\cal U}_a$. It has been proved so far only
in cases of repulsive interactions (e.g. $\kappa > 0$ in (NLS)$_n$ or
(R3)$_n$) where no solitary waves occur and where ${\cal U}_a$ can be
taken as the set of solutions of the free Schr\"odinger equation, or for small data. We
shall say nothing more about that problem in this paper. \\

We conclude that brief introduction to scattering theory by explaining
the difference between short and long range situations. The short range
case is the case where the interaction in the equation or system at
hand decays sufficiently fast at infinity in space and/or time for the
set of solutions of the underlying free system consisting of the free
Schr\"odinger equation and possibly of the free wave of Klein-Gordon
equation to be adequate as the set ${\cal U}_a$ of asymptotic motions.
The long range case is the complementary case where that set is
inadequate and has to be replaced by a set of modified asymptotic
motions. The modification includes in general (but not always, e.g.
(KGS)$_2$) the introduction of a phase in the asymptotic Schr\"odinger
function, and possibly the addition of correcting terms to the
Schr\"odinger and to the Wave or Klein-Gordon function. With respect to
that classification, (LS)$_n$ and (R3)$_n$ with potential (\ref{1.1e})
are short range for $\gamma > 1$ and long range for $\gamma \leq 1$,
(NLS)$_n$ is short range for $p > 1+2/n$ and long range for $p \leq 1
+ 2/n$, (WS)$_3$, (KGS)$_2$ and (MS)$_3$ are borderline long range,
namely analogous to (LS)$_n$ and (R3)$_n$ with $\gamma = 1$, while
(WS)$_n$, (KGS)$_n$ are short range for larger values of $n$. Finally
(Z)$_3$ is short range, while (Z)$_2$ is in a mixed situation and
exhibits some difficulties typical of the long range case. In this
paper, we shall review some of the results available for the systems
(WS)$_3$, (MS)$_3$, (Z)$_3$, (Z)$_2$ and (KGS)$_2$, in that order.\par

The construction of the wave operators and more precisely the local
Cauchy problem at infinity in the long range cases of the previous
nonlinear equations and systems has been treated essentially by two
methods, of which we shall concentrate on the first one. The second one
will only be quoted below for completeness. The first method is
intrinsically restricted to the case of small Schr\"odinger data and to
the borderline long range case. It applies to (NLS)$_n$ with $p =
1+2/n$ and $1 \leq n \leq 3$, to (R3)$_n$ with $V$ given by
(\ref{1.1e}) with $\gamma = 1$ and $n \geq 2$, to (WS)$_3$, (MS)$_3$ and
(KGS)$_2$. It can also be applied to (Z)$_n$ for $n = 2,3$. That method
was initiated by Ozawa in the case of (NLS)$_1$ \cite{25r}. It was
extended by Ozawa and one of the authors to (NLS)$_{2,3}$ and to (R3)$_n$
\cite{5r}, by Ozawa and Tsutsumi to (KGS)$_2$ and (Z)$_3$ \cite{27r}
\cite{28r} and by Tsutsumi to (MS)$_3$ \cite{36r}. In those early works
on (KGS)$_2$, (MS)$_3$ and (Z)$_3$ (in the last case for large data),
it was assumed in addition that the Fourier transform $\widehat{u}_+$ of the
Schr\"odinger asymptotic state $u_+$ satisfied a suitable support
condition, in order to cope with a difficulty coming from the
difference of asymptotic properties of the solutions of the
Schr\"odinger and of the wave or Klein-Gordon equations. Furthermore a
smallness assumption was made on the Klein-Gordon or Maxwell data in
the (KGS)$_2$ and (MS)$_3$ cases respectively. Recently, the problem
was revived in a series of papers by Shimomura \cite{29r}-\cite{34r}
who considered the (KGS)$_2$, (WS)$_3$, (MS)$_3$ and (Z)$_3$ systems.
The main achievement of that work is to eliminate the previous support
condition on $\widehat{u}_+$. This is made possible by using an improved form of
the asymptotic $u_a$ including an additional correction term which
partly cancels the free part of $A$ in the Schr\"odinger equation.
Furthermore in the (KGS)$_2$ case, the smallness condition of the KG
data is eliminated by going to (smaller) more regular function spaces
and by using a detailed asymptotic form for $(u_a, A_a)$ \cite{30r}
\cite{31r}. Finally the problem was revisited by the authors in the case of the (WS)$_3$, (MS)$_3$ and
(Z)$_{2,3}$ systems \cite{14r}-\cite{16r}. In that work, the smallness condition on the wave
or Maxwell field is eliminated for (WS)$_3$ and (MS)$_3$ respectively.
Furthermore, a refinement of the estimates and a more systematic use of
Strichartz inequalities allows for the use of (larger) less regular
function spaces. It is then possible to accomodate asymptotic $(u_a,
A_a)$ of low accuracy, so that the problem can be treated without the
Shimomura improved asymptotics and with much weaker assumptions on the
asymptotic state. \par

For completeness and although we shall not elaborate further on this
point in the present paper, we mention that the same problem, namely
the local Cauchy problem at infinity, can also be treated, at least for
the (R3)$_n$ equation and for the (WS)$_3$ and (MS)$_3$ systems, by a
more complex method where one first applies a phase-amplitude separation
to the Schr\"odinger function, inspired by previous work  by Hayashi
and Naumkin on the (R3)$_n$ equation \cite{18r} \cite{19r}. The main
interest of that method is to eliminate the smallness condition on the
Schr\"odinger function and to go beyond the $\gamma = 1$ borderline
case for the (R3)$_n$ equation. It has been applied by the
authors to the (R3)$_n$ equation \cite{8r} \cite{9r} with improvements
by Nakanishi eliminating a loss of regularity between the asymptotic
state and the solution \cite{23r} \cite{24r}, to the (WS)$_3$ system
\cite{10r} \cite{11r} and to the (MS)$_3$ system in a special case
\cite{12r}. \par

In this paper, as mentioned earlier, we shall review the main results
available by the first method for the systems (WS)$_3$, (MS)$_3$,
(Z)$_3$, (Z)$_2$ and (KGS)$_2$. Before describing the contents of this paper in
more detail however, we need to introduce some notation and to give
the principle of that method. The construction of solutions $(u, A)$ of
the system at hand with prescribed asymptotics $(u_a, A_a)$ is
performed in two steps. \\

\noi \underbar{Step 1.} One performs a change of variables, namely one
looks for $(u,A)$ in the form $(u, A) = (u_a + v, A_a + B)$ and one
studies the system satisfied by $(v, B)$. For instance in the (WS)$_3$
case, that auxiliary system takes the form  
\beq \label{1.6e}
\left \{ \begin{array}{l}  i\partial_t v = - (1/2)
\Delta v + A v + Bu_a - R_1 \\ \\ \sq B = - (|v|^2 + 2 \ {\rm Re} \ \overline{u}_a v) - R_2   \end{array} \right . 
 \eeq

\noi where the remainders $R_1$, $R_2$ are defined by
\beq \label{1.7e}
\left \{ \begin{array}{l}  R_1= i\partial_t u_a + (1/2)
\Delta u_a - A_a u_a \\ \\ R_2 = \sq A_a + |u_a|^2   \ . \end{array} \right . 
 \eeq

\noi The remainders measure the failure of the asymptotic form $(u_a,
A_a)$ to satisfy the original system. In particular the decay in time of
the remainders measures the quality of $(u_a, A_a)$ as an asymptotic
form. It plays an essential role in the problem. The first step of the
method consists in solving the auxiliary system for $(v, B)$ with $(v,
B)$ tending to zero at infinity in suitable norms under assumptions on
$(u_a, A_a)$ of a general nature, the most important of which being
decay assumptions on the remainders $R_1$ and $R_2$. That can be done
as follows. One first linearizes partly the auxiliary system satisfied
by $(v, B)$. For instance in the (WS)$_3$ case, the linearized
auxiliary system takes the form
 \beq \label{1.8e}
\left \{ \begin{array}{l}  i\partial_t v' = - (1/2)
\Delta v' + A v' + Bu_a - R_1 \\ \\ \sq B' = - (|v|^2 + 2 \ {\rm Re} \ \overline{u}_a v) - R_2  \ . \end{array} \right . 
 \eeq

\noi One solves the linearized system for $(v', B')$ with $(v', B')$
tending to zero at infinity, for fixed $(v, B)$ tending to zero at
infinity. This defines a map $\phi : (v, B) \to (v', B')$. One then
shows by a contraction method that the map $\phi$ has a fixed point in
some Banach space $X(I)$ of pairs $(v, B)$ of functions defined in $I =
[T, \infty )$ for $T$ sufficiently large. The definition of the space
$X(I)$ must include two kinds of ingredients, to be taylored after the
problem at hand, namely~:\par

(i) local space or space time regularity of the functions, in order to
cope with the nonlinearity of the relevant system. We shall always
assume that $v \in {\cal C}(I, H^k)$ and possibly include other norms
of equivalent homogeneity. As regards $B$, we shall include norms of
the same homogeneity as $(B, \partial_t B) \in {\cal C} (I,
\dot{H}^{\ell} \oplus \dot{H}^{\ell - 1})$ for some suitably high
$\ell$ and possibly for lower values. For such a choice of $X(I)$, we
shall say that $v$ is at level $k$, that $B$ is at level $\ell$ and
that $X(I)$ and the corresponding theory are at level $(k, \ell )$. \par

(ii) time decay of $(v, B)$ as $t \to \infty$ in the relevant norms.
Intuitively that time decay will be expressed by an exponent $\lambda$
so that the relevant norms decay essentially as $t^{-\lambda}$ (see
Section~2 below for a more precise formulation). In order to perform
Step 1, $\lambda$ will have to be sufficiently large, namely we shall
need a lower bound on $\lambda$. The theories described in this paper
will then be characterized in first rough approximation by their level
$(k, \ell )$ and their decay exponent $\lambda$. \\

\noi \underbar{Step 2.} That step obviously consists in constructing
asymptotic $(u_a, A_a)$ satisfying the assumptions needed for Step 1
and in particular the time decay assumptions of the remainders $R_1$
and $R_2$. That decay has to be sufficient to allow for the exponent
$\lambda$ needed for Step 1. It is therefore of interest to have the
lower bound in Step 1 as low as possible in order to accomodate the
most general (and therefore least accurate) possible $(u_a, A_a)$.\\

We can now describe the contents of this paper. In Section 2, we
consider the (WS)$_3$ system which is both the simplest one and the most
representative. After some preliminaries of general interest, which
will be used again in the subsequent sections, we implement Step 1
above at the lowest available level, namely $(k, \ell) = (0,1/2)$ with
lower bound $\lambda > 3/8$ (Proposition 2.1). We then describe the
simplest adequate choice of $(u_a, A_a)$ and the associated final
result, which has level $(k, \ell ) = (0, 1/2)$ and $\lambda = 1/2$
(Proposition 2.2). We then describe the Shimomura improved asymptotic
$(u_a, A_a)$ and the associated final result, which has again level
$(k, \ell ) = (0, 1/2)$ but now with $\lambda = 1$ up to logarithms
(Proposition 2.3). We then turn to more regular theories and in
particular we implement Step 1 at the level $(k, \ell ) = (2,1)$.
Remarkably enough the lower bound on $\lambda$ remains unchanged at
$\lambda > 3/8$ (Proposition 2.4). We then describe the final results
at that level, first with the simplest asymptotics and with $\lambda =
1/2$ (Proposition 2.5) and then with the improved asymptotics and
$\lambda = 1$ up to logarithms (Proposition 2.6). In Section 3, we
consider the (MS)$_3$ system in the Coulomb gauge. It turns out that
the lowest level available theory for that system is very similar to
although much more complicated than the level $(2,1)$ theory for
(WS)$_3$. We first implement Step 1 for (MS)$_3$ at the level $(k, \ell
) = (2, 3/2)$ with the same lower bound $\lambda > 3/8$ as for (WS)$_3$
(Proposition 3.1). We then describe the simplest adequate $(u_a, A_a)$
and the associated final result at level $(k, \ell ) = (2, 3/2)$ with
$\lambda = 1$ up to logarithms (Proposition 3.2). There is no need in
that case to use an improved asymptotic $(u_a, A_a)$ for the (MS)$_3$
system in the Coulomb gauge. We then comment briefly on the possible
extension of the result to other gauges and in particular to the
Lorentz gauge. In Section 4 we consider the (Z)$_n$ system, first for
$n = 3$. That system is short range and consequently no smallness
condition is needed on any of the data. We first implement Step 1 at
the lowest available level, namely $(k, \ell ) = (2, 1)$, with the
lower bound $\lambda > 1/4$ (Proposition 4.1). The simplest adequate
$(u_a, A_a)$ is then provided by the solutions of the free
Schr\"odinger and wave equations, and yields immediately the
corresponding final result with $\lambda = 1/2$ (Proposition 4.2),
while the improved asymptotic $(u_a, A_a)$ yields the corresponding
final result with $\lambda = 3/2$ (Proposition 4.3). We then turn to
the (Z)$_2$ system. We first implement Step 1 at the level $(k, \ell )
= (2,1)$ with $\lambda > 1/2$ (Proposition 4.4). That result however
requires again small Schr\"odinger data and the conditions on $A_a$
imply that the asymptotic state for $A$ has to be zero. We then state
the final result with $u_a$ a solution of the free Schr\"odinger
equation and $A_a = 0$, at the level $(k , \ell ) = (2,1)$ and with
$\lambda = 1$ (Proposition 4.5). In Section 5 we consider the
(KGS)$_2$ system. We implement Step~1 first at the lowest available
level, namely $(k, \ell ) = (0,1)$ with lower bound $\lambda > 1/2$
(Proposition 5.1) and then at the level $(k, \ell ) = (2, 1)$, which is
the same as that used for the (Z)$_2$ theory, with the same lower bound
$\lambda > 1/2$. We then describe the construction of the appropriate
$(u_a, A_a)$ performed in \cite{29r}, and we state the final results
in a qualitative way only, since additional work would be needed to
optimize the combination of that construction with the preceding
treatment of Step 1.\par

We emphasize again the fact that in the results described above, the
Schr\"odinger data have to satisfy a smallness condition, except in the
case of (Z)$_3$, that the wave, Maxwell and Klein-Gordon data can be
arbitrarily large except in the case of (Z)$_2$ where the asymptotic
state for $A$ has to be zero, and that no support assumption is made on
$\widehat{u}_+$. \par

We conclude this introduction by giving some notation which will be
used freely in this paper. We denote by $F$ the Fourier transform in ${I\hskip-1truemm R}^n$, with $Ff = \widehat{f}$ for any function $f$, by
$\parallel \cdot \parallel_r$ the norm in $L^r = L^r({I\hskip-1truemm R}^n)$, $1 \leq r
\leq \infty$, $n =2,3$, and by $<\cdot , \cdot >$ the scalar product in
$L^2$. Beyond the standard Sobolev spaces $H^k$ defined for any $k \in
{I\hskip-1truemm R}$ by 
$$H^k = \left \{ u :\ \parallel u; H^k\parallel \ = \ \parallel <\omega >^k u \parallel_2\ < \infty \right \}$$

\noi where $\omega = (- \Delta)^{1/2}$ and $< \cdot > = (1 + |\cdot
|^2)^{1/2}$ and their homogeneous versions
$$\dot{H}^k = \left \{ u :\ \parallel u; \dot{H}^k\parallel \ = \ \parallel \omega^k u \parallel_2\ < \infty \right \}\ ,$$

\noi we shall use the spaces $H^{k,s}$ defined for $k, s \in
{I\hskip-1truemm R}$ by 
$$H^{k,s} = \left \{ u :\ \parallel u; H^{k,s}\parallel \ = \ \parallel <x>^s <\omega >^k u \parallel_2\ < \infty \right \}$$

\noi so that $H^{k,0} = H^k$ and $H^{0,k} = FH^k$. We shall also use
the Sobolev spaces $W_r^k$ defined for $1 \leq r \leq \infty$ and for
$k$ a non negative integer by 
$$W_r^k = \left \{ u : \parallel u; W_r^k\parallel\ = \sum_{\alpha : 0 \leq |\alpha | \leq k} \parallel \partial_x^{\alpha} u \parallel_r \ < \infty \right \}$$

\noi so that $W_2^k = H^k$. For any interval $I$, for any Banach space
$X$ and for any $q$, $1 \leq q \leq \infty$, we denote by $L^q(I, X)$
(resp. $L_{loc}^q(I, X)$) the space of $L^q$ integrable (resp. locally
$L^q$ integrable) functions from $I$ to $X$ if $q < \infty$ and the
space of measurable essentially bounded (resp. locally essentially
bounded) functions from $I$ to $X$ if $q = \infty$.

\mysection{The Wave-Schr\"odinger system (WS)$_3$}
\hspace*{\parindent} In this section we review the main results available on the local Cauchy problem at infinity for the (WS)$_3$ system
\beq \label{2.1e}
\left \{ \begin{array}{l}  i\partial_t u = - (1/2)
\Delta u + A u \\ \\ \sq A = - |u|^2  \ . \end{array} \right . 
 \eeq

That system is known to be globally well posed in the energy space
\cite{1r} \cite{4r}. \par

The exposition is based mostly on \cite{14r} but includes also some
results of \cite{32r} rephrased in the framework of \cite{14r}. We
follow the sketch given in the introduction and we first consider Step
1. For a given asymptotic $(u_a, A_a)$, we look for $(u, A)$ in the
form $(u, A) = (u_a + v, A_a + B)$ and we try to solve the auxiliary
system (\ref{1.6e}) for $(v, B)$ with $(v, B)$ tending to zero at
infinity. For that purpose we first solve the partly linearized system
(\ref{1.8e}) for $(v', B')$ and we try to prove that the map $\phi :
(v, B) \to (v', B')$ thereby defined is a contraction in a suitable
space $X(I)$ with $I = [T, \infty )$. The crucial point of Step 1 is
the choice of $X(I)$ and that choice is dictated by the available
estimates. We use three types of estimates. \par

(i) \underbar{$L^2$ or energy estimates.} If $v$ satisfies the Schr\"odinger equation 
\beq
\label{2.2e}
i\partial_t v = - (1/2)\Delta v + A v + f
\eeq

\noi then
$$\partial_t \parallel v \parallel_2^2\ = 2\ {\rm Im} \ <v,f>$$

\noi so that if $v$ tends to zero at infinity, then
\beq
\label{2.3e}
\parallel v(t) \parallel_2\ \leq \ \parallel f;L^1([t, \infty ), L^2)\parallel \ .
\eeq

\noi An important feature of that estimate is that it does not involve
$A$. This is the key to the elimination of smallness conditions on $A$
in the theory. More generally, if $\partial = \nabla$ or $\partial_t$,
from the equation
\beq
\label{2.4e}
i\partial_t \partial v = - (1/2)\Delta \partial v + A \partial v + v\partial A + \partial f
\eeq 

\noi we obtain 
\beq
\label{2.5e}
\parallel \partial v(t) \parallel_2\ \leq \ \parallel v \partial A + \partial f;L^1([t, \infty ), L^2)\parallel 
\eeq

\noi if $\partial v$ tends to zero at infinity. The estimate
(\ref{2.3e}) has level $k = 0$, while (\ref{2.5e}) has level $k = 1$
if $\partial = \nabla$ and level $k = 2$ if $\partial = \partial_t$.
\par

We shall also use the energy estimate for the wave equation. If $B$
satisfies 
\beq
\label{2.6e}
\sq B = g
\eeq

\noi and tends to zero at infinity, then
\beq
\label{2.7e}
\parallel \partial_t B(t) \parallel_2\ \vee \ \parallel \nabla B(t) \parallel_2\ \leq \ \parallel g;L^1([t, \infty ), L^2)\parallel \ .
\eeq

\noi That estimate has level $\ell = 1$.\par

(ii) \underbar{Strichartz inequalities for the Schr\"odinger equation.}\par

We recall those inequalities for completeness, in dimension $n\geq 2$ (see \cite{3r} \cite{20r} \cite{37r} and references). A pair of exponents $q$, $r$ with $2 \leq q$, $r \leq \infty$ is called admissible if
\beq
\label{2.8e}
\begin{array}{ll} 0 \leq 2/q = n/2 - n/r &\leq 1 \quad \hbox{for $n \geq 3$}\\
&\\
&<1 \quad \hbox{for $n = 2$} \ . \\
\end{array}
\eeq

\noi {\bf Lemma 2.1.} {\it Let $(q_i, r_i)$, $i = 1,2$, be two admissible pairs. Let $v$ satisfy the equation
$$i \partial_t v = - (1/2) \Delta v + f$$

\noi in some interval $I$ with $v(t_0) = v_0$ for some $t_0 \in I$. Then the following estimates hold~:
\beq
\label{2.9e}
\parallel v; L^{q_1}(I, L^{r_1})\parallel \ \leq \ C \left ( \parallel v_0\parallel_2 \ + \ \parallel f; L^{\overline{q}_2}(I, L^{\overline{r}_2})\parallel \right )
\eeq

\noi where $C$ is a constant independent of $I$, and with $1/p +
1/\overline{p} = 1$. }\\

\noi The Strichartz inequalities do not involve any derivative and are
written in (\ref{2.9e}) at the level $k=0$. Their extension to general
level $k$ is obvious.\par

(iii ) \underbar{Strichartz inequalities for the wave equation.}\par

We shall use and therefore we state those inequalities only in
dimension 3 and for a special case of exponents (see \cite{7r}
\cite{20r} and references for the general case).\\

\noi {\bf Lemma 2.2.} {\it Let $n=3$ and let $B$ satisfy the equation
(\ref{2.6e}) in some interval $I$ with $B(t_0) = B_0$, $\partial_t
B(t_0) = B_1$ for some $t_0 \in I$. Then the following estimates hold~:
\beq
\label{2.10e}
\parallel B; L^4(I, L^4)\parallel \ \leq \ C \left ( \parallel \omega^{1/2} B_0 \parallel_2\ + \ \parallel \omega^{-1/2}B_1 \parallel_2 \ + \ \parallel g; L^{4/3}(I, L^{4/3})\parallel \right ) \ ,
\eeq
\bea
\label{2.11e}
&&\parallel \nabla B; L^4(I, L^4)\parallel \ \vee \ \parallel \partial_t B;L^4(I,L^4)\parallel \nn \\
&&\leq \ C \left ( \parallel \omega^{3/2} B_0 \parallel_2\ + \ \parallel \omega^{1/2}B_1 \parallel_2 \ + \ \parallel \nabla g; L^{4/3}(I, L^{4/3})\parallel \right ) 
\eea

\noi where $C$ is a constant independent of $I$, and $\omega = (-
\Delta )^{1/2}$.}\\

The estimates (\ref{2.10e}) and (\ref{2.11e}) have level $\ell = 1/2$
and $\ell = 3/2$ respectively. \par

The definition of the space $X(I)$ will involve a suitable time decay
which will be characterized by a function $h \in {\cal C}([1, \infty ),
{I\hskip-1truemm R}^+)$ such that for a suitable $\lambda > 0$, the
function $\overline{h}(t) = t^{\lambda}h(t)$ is nonincreasing and tends
to zero at infinity. This means in particular that $h$ decreases
slightly faster than $t^{-\lambda}$.\par

We can now state the result concerning Step 1 for the (WS)$_3$ system
at the lowest available level, which is $(k, \ell ) = (0, 1/2)$. That
result is the prototype of all the results concerning Step 1 contained
in this paper. The relevant function space is defined by
\bea
\label{2.12e}
&&X(I) = \Big \{ (v, B):v\in {\cal C}(I, L^2), \parallel (v, B);X(I)\parallel \ \equiv \ \mathrel{\mathop {\rm Sup}_{t \in I }}\ h(t)^{-1} \nn \\
&&\left ( \parallel v(t)\parallel_2 \ + \ \parallel v; L^{8/3}(J,L^4)\parallel\ + \ \parallel B; L^4(J, L^4)\parallel \right ) < \infty \Big \} 
\eea

\noi for any interval $I \subset [1, \infty )$, with $J = [t, \infty )
\cap I$. The pair of exponents $(8/3, 4)$ is Schr\"odinger admissible, and
the norm involves the $L^2$ norm and a Strichartz norm of level $k=0$
for $v$, and a Strichartz norm of level $\ell = 1/2$ for $B$. The
result is Proposition 2.2 of \cite{14r} and can be stated as follows.
\\

\noi {\bf Proposition 2.1.} {\it Let $h$ be defined as above with
$\lambda = 3/8$ and let $X(\cdot )$ be defined by (\ref{2.12e}). Let
$(u_a, A_a)$ be sufficiently regular (for the following estimates to
make sense) and satisfy the estimates
\beq
\label{2.13e}
\parallel u_a(t) \parallel_4\ \leq c_4 \ t^{-3/4}\ ,
\eeq
\beq
\label{2.14e}
\parallel A_a(t) \parallel_{\infty} \ \leq a\ t^{-1}\ ,
\eeq
\beq
\label{2.15e}
\parallel R_1;L^1([t, \infty ), L^2) \parallel\ \leq r_1\ h(t)\ ,
\eeq
\beq
\label{2.16e}
\parallel R_2;L^{4/3}([t, \infty ), L^{4/3}) \parallel\ \leq r_2\ h(t)\ ,
\eeq

\noi for some constants $c_4$, $a$, $r_1$ and $r_2$ with $c_4$
sufficiently small and for all $t \geq 1$. Then there exists $T$, $1
\leq T < \infty$ and there exists a unique solution $(v, B)$ of the
system (\ref{1.6e}) in the space $X([T, \infty ))$.} \\

\noi {\bf Remark 2.1.} The time decay assumed in (\ref{2.13e})
(\ref{2.14e}) for $(u_a, A_a)$ is the optimal time decay that can be
obtained for solutions of the free Schr\"odinger and wave equations
respectively. We shall always make assumptions of this type in all
subsequent results concerning Step 1.\\

\noi {\bf Remark 2.2.} There is an absolute smallness condition on
$u_a$ through $c_4$ but none other. In particular $A_a$ can be
arbitrarily large.\\

\noi {\bf Sketch of proof.} We follow the sketch given in the
introduction. Let $1 \leq T \leq \infty$ and let $(v, B) \in X([T,
\infty))$, so that $(v, B)$ satisfies
\beq
\label{2.17e}
\left \{ \begin{array}{l} \parallel v(t) \parallel_2 \ \leq N_0\ h(t)\\
\\
\parallel v ; L^{8/3}([t, \infty), L^4) \parallel \ \leq N_1\ h(t)\\
\\
\parallel B ; L^{4}([t, \infty), L^4) \parallel \ \leq N_2\ h(t)
\end{array}\right .\eeq
 
\noi for some constants $N_i$ and for all $t \geq T$. We first solve
the linearized system (\ref{1.8e}) for $(v',B')$ in $X([T, \infty))$.
That can be done by first solving that system in $X([T, t_0])$ with
initial condition $(v', B')(t_0) = 0$ for some large $t_0 > T$, which
is an easy (linear) Cauchy problem with finite initial time, and then
taking the limit of the solution thereby obtained when $t_0 \to
\infty$. An essential point for taking that limit is to derive
estimates of the solution in $X(I)$ for $I = [T, t_0]$ that are uniform
in $t_0$. We define 
\beq
\label{2.18e}
\left \{ \begin{array}{l}N'_0 =  \ \displaystyle{\mathrel{\mathop {\rm Sup}_{t\in I}}}\ h(t)^{-1} \parallel v'(t) \parallel_2
\\ 
\\
N'_1 =  \ \displaystyle{\mathrel{\mathop {\rm Sup}_{t\in I}}}\ h(t)^{-1} \parallel v';L^{8/3}(J,L^4) \parallel
\\
\\
N'_2 =  \ \displaystyle{\mathrel{\mathop {\rm Sup}_{t\in I}}}\ h(t)^{-1} \parallel B';L^{4}(J,L^4) \parallel\\
\end{array} \right .
\eeq

\noi where $J = [t, \infty ) \cap I$. Using an $L^2$ estimate of the
type (\ref{2.3e}) and Strichartz estimates of the type (\ref{2.9e})
(\ref{2.10e}), one can show that  
\beq
\label{2.19e}
\left \{ \begin{array}{l} N'_0 \leq C_0 \left ( c_4\ N_2 + r_1 \right ) \\
\\
N'_1 \leq C_1 \left ( c_4 \ N_2 + r_1\right ) \left ( 1 + a + N_2\ \overline{h}(T)\right ) \\
\\
N'_2 \leq C_2 \left (  c_4\ N_0 + r_2 + N_0\ N_1 \ \overline{h}(T) \right )  \ ,
\end{array} \right .
\eeq

\noi where the $C_i$, $0 \leq i \leq 2$, are absolute constants. Using
the fact that those estimates are uniform in $t_0$, one can easily take
the limit $t_0 \to \infty$ of the previous solution, thereby obtaining
a solution $(v', B') \in X([T, \infty))$ also satisfying the estimates
(\ref{2.19e}) with the $N'_i$ defined by (\ref{2.18e}), now with $I =
[T, \infty )$. At this stage we have completed the construction of the
map $\phi : (v, B) \to (v', B')$ and we have proved in addition that
this map is bounded in $X(I)$. We next show that $\phi$ is a
contraction on a suitable closed subset ${\cal R}$ of $X(I)$ for $T$
sufficiently large. We define ${\cal R}$ by (\ref{2.17e}).\par

From (\ref{2.19e}), it follows that ${\cal R}$ is stable under $\phi$
provided
\beq
\label{2.20e}
\left \{ \begin{array}{l} 
C_0 \left ( c_4 \ N_2 + r_1 \right ) \leq N_0 \\ \\ C_1 \left ( c_4\ N_2 + r_1 \right ) \left ( 1 + a + N_2 \ \overline{h}(T)\right ) \leq N_1 \\ \\ C_2 \left ( c_4\ N_0 + r_2 + N_0 \ N_1\ \overline{h}(T)\right ) \leq N_2 
\end{array} \right .
\eeq

\noi which can be ensured under the smallness condition $C_0C_2c_4^2 < 1$ by choosing the $N_i$ according to 
\beq
\label{2.21e}
\left \{ \begin{array}{l} 
N_0 = C_0 \left ( c_4 \ N_2 + r_1 \right ) \\ \\ N_1 =  C_1 \left ( c_4\ N_2 + r_1 \right ) (2+a)  \\ \\ N_2 =  C_2 \left ( c_4\ N_0 + r_2 + 1 \right )  
\end{array} \right .
\eeq

\noi and by taking $T$ sufficiently large so that
\beq
\label{2.22e}
N_2\ \overline{h} (T) \leq 1 \ , \ N_0\ N_1\ \overline{h} (T) \leq 1 \ .
\eeq

It remains to prove that $\phi$ is a contraction on ${\cal R}$. This is
done by estimating the difference of two solutions $(v', B')$ of
(\ref{1.8e}) corresponding to two different $(v, B)$. The estimates are
minor variants of (\ref{2.19e}) (see the proof of Proposition 2.2 in
\cite{14r} for details). \par \nobreak \hfill $\sq$\par

We now turn to Step 2, namely to the construction of $(u_a, A_a)$
satisfying the assumptions of Proposition 2.1. One sees immediately
that taking for $(u_a, A_a)$ a pair of solutions $(u_0, A_0)$ of the
free Schr\"odinger and wave equations is inadequate, since in that case
$R_2 = |u_0|^2$ so that at best 
$$\parallel R_2\parallel_{4/3} \ \leq \ \parallel u_0 \parallel_2\ \parallel u_0 \parallel_4\ \leq C\ t^{-3/4} \notin L^{4/3} \ .$$

However, with the weak time decay allowed by Proposition 2.1, the
simplest modification available in the literature suffices. We describe
that choice first. We choose
\beq
\label{2.23e}
A_a = A_0 + A_1
\eeq
\beq
\label{2.24e}
A_0 = \cos \omega t \ A_+ + \omega^{-1} \sin \omega t \ \dot{A}_+
\eeq
\beq
\label{2.25e}
A_1= \int_t^{\infty} dt' \omega^{-1} \sin (\omega (t-t'))|u_a (t')|^2 
\eeq

\noi where $\omega = (-\Delta )^{1/2}$ and $(A_+, \dot{A}_+)$ is the
asymptotic state of $A$. That choice ensures that $\sq A_0 = 0$,
$\sq A_a = \sq A_1 = - |u_a|^2$ so that $R_2 = 0$. We next choose
$u_a$. The Schr\"odinger group (\ref{1.5e}) admits the well known
factorisation
\beq
\label{2.26e}
U(t) \equiv \exp (i(t/2)\Delta ) = M\ D\ F\ M
\eeq

\noi where
\beq
\label{2.27e}
M \equiv M(t) = \exp (ix^2/2t) 
\eeq
\beq
\label{2.28e} 
D(t) = (it)^{-3/2} D_0(t) \quad , \quad ( D_0 (t) f) (x) = f(x/t) 
\eeq

\noi and $F$ is the Fourier transform. We choose
\beq
\label{2.29e}
u_a = MD \exp (- i \varphi ) \widehat{u}_+
\eeq 

\noi where $\varphi$ is a real phase, $\widehat{u}_+ = Fu_+$ and $u_+$ is the
asymptotic state of $u$. Substituting (\ref{2.29e}) into (\ref{2.25e})
yields
\beq
\label{2.30e}
A_1(t) = t^{-1} D_0(t) \widetilde{A}_1\ ,
\eeq

\noi where
\beq
\label{2.31e}
\widetilde{A}_1 = - \int_1^{\infty} d\nu \ \nu^{-3} \omega^{-1} \sin (\omega (\nu - 1)) D_0 (\nu ) |\widehat{u}_+|^2 \ .
\eeq

\noi In particular $\widetilde{A}_1$ is constant in time. Substituting
(\ref{2.23e}) (\ref{2.29e}) (\ref{2.30e}) into the definition
(\ref{1.7e}) of $R_1$ and using the commutation relation
\beq
\label{2.32e}
\left ( i \partial_t + (1/2) \Delta \right ) MD = MD \left ( i \partial_t + (2t^2)^{-1}\Delta \right ) 
\eeq

\noi we obtain
\beq
\label{2.33e}
R_1 = \widetilde{R}_1 - A_0\ u_a
\eeq

\noi where
\bea
\label{2.34e}
\widetilde{R}_1 &=&\left ( i \partial_t + (1/2) \Delta - A_1 \right ) MD \exp (-i \varphi ) \widehat{u}_+ \nn \\
&=& MD \left ( i \partial_t + (2t^2)^{-1} \Delta - t^{-1} \widetilde{A}_1 \right ) \exp (- i \varphi ) \widehat{u}_+ \ .
\eea

\noi We finally use $\varphi$ to cancel the long range term $\widetilde{A}_1$ by taking 
\beq
\label{2.35e}
\varphi = (\ell n\ t) \widetilde{A}_1 
\eeq

\noi thereby obtaining
\beq
\label{2.36e}
\widetilde{R}_1 = (2t^2)^{-1} MD\Delta \exp (-i \varphi ) \widehat{u}_+ \ .
\eeq

In order to show that $(u_a, A_a)$ satisfies the assumptions needed for
Step 1, we need the well known dispersive estimate for the
Schr\"odinger equation, namely
\beq
\label{2.37e}
\parallel U(t) u_+ \parallel_r \ \leq \ (2 \pi |t|)^{- \delta (r)} \parallel u_+ \parallel_{\overline{r}}
\eeq

\noi where $2 \leq r \leq \infty$, $\delta (r) = n/2 - n/r$ and $1/r +
1/\overline{r} = 1$, and some general estimates of solutions of the
free wave equation, which we recall for completeness \cite{35r}. \\

\noi {\bf Lemma 2.3.} {\it Let $A_0$ be defined by (\ref{2.24e}) and
let $\ell \geq 0$ be an integer. Let $(A_+, \dot{A}_+)$ satisfy the
conditions
\beq
\label{2.38e}
A_+, \omega^{-1} \dot{A}_+ \in H^{\ell} \ .
\eeq

\noi Then $A_0$ satisfies the estimates
\beq
\label{2.39e}
\left \{ \begin{array}{l} \parallel A_0(t) ; H^{\ell} \parallel \ \leq C \ , \\ \\ \parallel \partial_t A_0(t) ; H^{\ell - 1} \parallel \ \leq C \quad \hbox{\it for $\ell \geq 1$} \ .
\end{array} \right .
\eeq

\noi If in addition $(A_+, \dot{A}_+)$ satisfies the conditions
\beq
\label{2.40e}
\nabla^2 A_+\ , \ \nabla \dot{A}_+ \in W_1^{\ell} \ ,
\eeq

\noi then $A_0$ satisfies the estimates
\beq
\label{2.41e}
\left \{ \begin{array}{l} \parallel A_0(t) ; W_r^{\ell} \parallel \ \leq C \ t^{-1 + 2/r}\ , \\ \\ \parallel \partial_t A_0(t) ; W_r^{\ell - 1} \parallel \ \leq C \ t^{-1 + 2/r} \quad \hbox{\it for $\ell \geq 1$} 
\end{array} \right .
\eeq

\noi for $2 \leq r \leq \infty$ and for all $t \in {I\hskip-1truemm R}$.} \\

Using Lemma 2.3 and Sobolev inequalities, one can easily derive the
estimates of $(u_a, A_a)$ needed in Proposition 2.1 (see Proposition
3.1, part (1) in \cite{14r}). \\

\noi {\bf Lemma 2.4.} {\it Let $u_+ \in H^{0,2}$ and let $(A_+,
\dot{A}_+)$ satisfy (\ref{2.38e}) (\ref{2.40e}) with $\ell = 0$. Then the
following estimates hold
\beq
\label{2.42e}
\parallel u_a(t)\parallel_r\ \leq t^{-\delta (r)} \parallel \widehat{u}_+ \parallel_r \qquad \hbox{\it for $2 \leq r \leq \infty$} \ ,
\eeq

\noi where $\delta (r) = 3/2 - 3/r$,}
\beq
\label{2.43e}
\parallel A_a(t)\parallel_{\infty} \ \leq a\ t^{-1} \ ,
\eeq
\beq
\label{2.44e}
\parallel \widetilde{R}_1(t)\parallel_2 \ \leq r_1\ t^{-2} (2 + \ell n \ t)^2\ ,
\eeq
\beq
\label{2.45e}
\parallel R_1(t)\parallel_2 \ \leq r_1\ t^{-3/2}\ .
\eeq

Putting together Proposition 2.1 and Lemma 2.4, we obtain the first
(and simplest) final result for the system (WS)$_3$ at the level $(k,
\ell ) = (0,1/2)$ (see Proposition 1.1, part (1) in \cite{14r}). \\

\noi {\bf Proposition 2.2.} {\it Let $h(t) = t^{-1/2}$ and let $X(\cdot
)$ be defined by (\ref{2.12e}). Let $(u_a, A_a)$ be defined by
(\ref{2.29e}) (\ref{2.35e}) (\ref{2.23e}) (\ref{2.24e}) (\ref{2.30e})
(\ref{2.31e}). Let $u_+ \in H^{0,2}$ with $c_4 = \break\noindent \parallel \widehat{u}_+
\parallel_4$ sufficiently small. Let $A_+$, $\omega^{-1} \dot{A}_+ \in
L^2$ and $\nabla^2A_+$, $\nabla \dot{A}_+ \in L^1$. Then there exists
$T$, $1 \leq T < \infty$, and there exists a unique solution $(u, A)$
of the (WS)$_3$ system (\ref{2.1e}) such that $(v, B) \equiv (u - u_a,
A- A_a) \in X([T, \infty ))$.}\\

It is apparent from Lemma 2.4, especially from (\ref{2.44e})
(\ref{2.45e}) that the time decay $t^{-1/2}$ in Proposition 2.2 comes
from the term $A_0u_a$ in $R_1$ (see (\ref{2.33e})) which is estimated
only as
$$\parallel A_0 u_a \parallel_2 \ \leq \ \parallel A_0 \parallel_2\ \parallel u_a \parallel_{\infty} \ \leq C\ t^{-3/2}$$

\noi whereas $\widetilde{R}_1$ would allow for the better decay
$t^{-1}(2 + \ell n \ t)^2$. In \cite{32r}, Shimomura has proposed an
improved form of $u_a$ and has obtained that improved decay on a
subspace of more regular and decaying asymptotic states. We now
describe the results that can be obtained thereby. We keep the same
$A_a$ as before, but we now take $u_a = u_1 + u_2$ where $u_1$ is the
previous choice of $u_a$, namely
\beq
\label{2.46e}
u_1 = MD \exp (- i \varphi ) \widehat{u}_+
\eeq

\noi and $u_2$ is chosen in such a way that $(i \partial_t + (1/2)
\Delta ) u_2$ in $R_1$ partly cancels the term $A_0 u_1$. For that
purpose, it is appropriate to look for $u_2$ in the form $u_2 = fu_1$
for a suitable real function $f$, so that now
\beq
\label{2.47e}
u_a = (1 + f) u_1
\eeq

\noi with $u_1$ defined by (\ref{2.46e}). Substituting (\ref{2.47e}) into the definition of $R_1$ yields
\bea
\label{2.48e}
R_1 &=& \left ( i \partial_t + (1/2) \Delta - A_a \right ) (1 + f) u_1\nn \\
&=& (1 + f) \left ( i \partial_t + (1/2) \Delta - A_1 \right ) u_1 + \left ( (1/2) \Delta f - A_0 \right ) u_1 \nn\\
&&- f A_0 u_1 + \nabla f \cdot \nabla u_1 + i \left ( \partial_t f\right ) u_1 \ .
\eea

\noi We now choose
\beq
\label{2.49e}
f = 2 \Delta^{-1} A_0 \ .
\eeq

\noi Making that choice and using the operators 
\beq
\label{2.50e}
J = x + it\nabla \quad , \quad P = t \partial_t + x \cdot \nabla \ ,
\eeq

\noi we obtain (see (\ref{2.34e}))
\beq
\label{2.51e}
R_1 = (1 + f) \widetilde{R}_1 - f A_0 u_1 - i t^{-1} \nabla f \cdot Ju_1 + i t^{-1} (Pf) u_1 \ .
\eeq

\noi Now $f$ and $\nabla f$ are solutions of the free wave equation. From the commutation rule
\beq
\label{2.52e}
\sq P = (P+2) \sq
\eeq

\noi it follows that the same is true for $Pf$. By Lemma 2.3, one can
ensure that $f$, $\nabla f$ and $Pf$ are uniformly bounded in $L^2$ by
making suitable assumptions on their initial data, which reduce to
assumptions on $(A_+, \dot{A}_+)$. On the other hand from the
commutation relation 
\beq
\label{2.53e}
JMD = iMD\nabla
\eeq

\noi it follows that
\beq
\label{2.54e}
\parallel Ju_1 \parallel_{\infty} \ \leq C\ t^{-3/2} (2 + \ell n\ t) \ .
\eeq

\noi The last three terms in (\ref{2.51e}) are then estimated by
\beq
\label{2.55e}
\left \{ \begin{array}{l}
\parallel fA_0 u_1 \parallel_2\ \leq \ \parallel f \parallel_2 \ \parallel A_0 \parallel_{\infty} \ \parallel u_1 \parallel_{\infty}\ \leq C\ t^{-5/2}\\ \\
t^{-1}\parallel \nabla f \cdot Ju_1 \parallel_2\ \leq \ t^{-1} \parallel \nabla f \parallel_2 \ \parallel Ju_1\parallel_{\infty}\ \leq \ C \ t^{-5/2} (2 + \ell n\ t)\\ \\
t^{-1}\parallel (P f ) u_1 \parallel_2\ \leq \ t^{-1} \parallel P f \parallel_2 \ \parallel u_1 \parallel_{\infty} \ \leq \ C \ t^{-5/2}\\
\end{array} \right .
\eeq

\noi and would therefore allow for a decay $h(t) = t^{-3/2} (2 + \ell
n\ t)$, so that now the dominant term in $R_1$ is $\widetilde{R}_1$,
thereby allowing for a decay $h(t) = t^{-1} (2 + \ell n\ t)^2$.\par

With the new choice (\ref{2.46e}) (\ref{2.47e}) of $u_a$ and with the
same choice of $A_a$ as before, the remainder $R_2$ is no longer zero
and becomes
\beq
\label{2.56e}
R_2 = \left ( f^2 + 2 f\right ) |u_1|^2 \ .
\eeq

\noi This can be estimated as
$$\parallel R_2\parallel_{4/3} \ \leq \ \left ( \parallel f \parallel_{\infty} + 2 \right ) \parallel f\parallel_2 \ \parallel u_1 \parallel_{\infty}\ \parallel u_1\parallel_4 \ \leq C\ t^{-3/2-3/4}$$

\noi so that
\beq
\label{2.57e}
\parallel R_2;L^{4/3} ([t, \infty ), L^{4/3}) \parallel\ \leq C\ t^{-3/2} \ ,
\eeq

\noi which would allow for a decay $h(t) = t^{-3/2}$. \par

Making the appropriate assumptions required to implement the previous
estimates, one obtains the following final result with improved decay.\\

\noi {\bf Proposition 2.3.} {\it Let $h(t) = t^{-1}(2 + \ell n\ t)^2$ and let $X(\cdot
)$ be defined by (\ref{2.12e}). Let $(u_a, A_a)$ be defined by
(\ref{2.46e}) (\ref{2.47e}) (\ref{2.49e}) and (\ref{2.35e}) (\ref{2.23e})
(\ref{2.24e}) (\ref{2.30e}) (\ref{2.31e}). Let $u_+ \in H^{0,2}$ with $x u_+ \in L^1$ and with $c_4 = \ \parallel \widehat{u}_+
\parallel_4$ sufficiently small. Let $(A_+, \dot{A}_+)$ satisfy
\bea
\label{2.58e}
&&A_+, \omega^{-1} \dot{A}_+ \in \dot{H}^{-2} \cap L^2\quad , \quad \nabla^2 A_+, \nabla \dot{A}_+ \in L^1\ ,\nn \\
&&x \cdot \nabla A_+, \omega^{-1} x \cdot \nabla \dot{A}_+ \in \dot{H}^{-2} \ .
\eea
\noi  Then :\par

(1) There exists
$T$, $1 \leq T < \infty$, and there exists a unique solution $(u, A)$
of the (WS)$_3$ system (\ref{2.1e}) such that $(v, B) \equiv (u - u_a,
A- A_a) \in X([T, \infty ))$.\par

(2) There exists
$T$, $1 \leq T < \infty$, and there exists a unique solution $(u, A)$
of the (WS)$_3$ system (\ref{2.1e}) such that $(u - u_1,
A- A_a) \in X([T, \infty ))$. One can take the same $T$ and the solution $(u,A)$ is the same as in Part (1).}\\

Part (2) of the Proposition follows immediately from Part (1) and from
the fact that $(fu_1, 0) \in X([T, \infty ))$ for any $T \geq 1$ under
the assumptions made on $(u_+, A_+, \dot{A}_+)$. Thus the correcting
term $u_2 = fu_1$ does not change the asymptotics. It plays however an
essential role in estimating the time derivative of $v$. The same
situation will appear repeatedly, in particular in Propositions 2.6 and
4.3 below.\par

The subspace of asymptotic states $(u_+, A_+,\dot{A}_+)$ giving rise to
the improved decay of Proposition 2.3 is specified by the additional
assumptions of Proposition 2.3 as compared with those of Proposition
2.2. They include the condition $xu_+ \in L^1$, which expresses a space
decay dimensionally stronger by one half power of $x$ than $u_+ \in
H^{0,2}$, the condition that $A_+$, $\omega^{-1} \dot{A}_+ \in
\dot{H}^{-2}$ and a similar condition on $x \cdot \nabla A_+$,
$\omega^{-1} x \cdot \nabla \dot{A}_+$. The condition on $A_+$,
$\omega^{-1}\dot{A}_+$ expresses a space decay dimensionally stronger
than the corresponding $L^2$ condition by two powers of $x$, but that
condition cannot be ensured by space decay alone. In fact
$$\omega^{-2} A_+ = (4 \pi |x|)^{-1} \ *\ A_+$$

\noi and this can never be in $L^2$ unless the space integral of $A_+$
vanishes. The situation is even worse for $\dot{A}_+$. As an
illustration, we now give conditions on $A_+$, $\dot{A}_+$ in terms of
space decay and vanishing of suitable moments which suffice to ensure
the $\dot{H}^{-2}$ conditions. Similar results apply to $x\cdot \nabla
A_+$, $x \cdot \nabla \dot{A}_+$. The following result is a special
case of Lemma 3.5 in \cite{11r}.\\

\noi {\bf Lemma 2.5.} {\it Let $A_+$, $\omega^{-1} \dot{A}_+ \in L^2$.
Assume in addition that
$$<x>^{1/2 + \varepsilon} A_+ \in L^1\ , \ <x>^{3/2 + \varepsilon} \dot{A}_+ \in L^1 \quad \hbox{\it for some $\varepsilon > 0$}\ , $$ 
$$ \int dx\ A_+ = \int dx \dot{A}_+ = \int dx\ x\ \dot{A}_+ = 0 \ .$$

\noi Then $A_+$, $\omega^{-1} \dot{A}_+ \in \dot{H}^{-2}$.}\\

\noi {\bf Remark 2.3.} In Proposition 2.3 as in Proposition 2.2, there
is again a discrepancy between the time decay allowed by
$\widetilde{R}_1$ and by the other terms of $R_1$, but now in the
opposite direction since now $\widetilde{R}_1$ allows for $t^{-1}(2 +
\ell n\ t)^2$ whereas the other terms allow for $t^{-3/2}(2 + \ell n\
t)$. This can be amended in two ways. On the one hand one can weaken
the assumptions on $A_+$, $\dot{A}_+$ so that the last three terms of
$R_1$ have no better decay than $\widetilde{R}_1$. On the other hand
one can replace $u_1$ and in particular $\varphi$ by better
approximations in order to improve the time decay of $\widetilde{R}_1$
to $t^{-3/2}$, possibly up to logarithms. We refer to Section 3 of
\cite{11r} where similar questions are treated in a slightly different
context.\\

We have considered the (WS)$_3$ system so far at the lowest convenient
level of regularity $(k, \ell ) = (0, 1/2)$. However there is no
difficulty in constructing theories similar to the previous one at
higher levels of regularity. It suffices to include suitable higher
norms in the definition of $X(\cdot )$ and to estimate those norms in
the proofs by straightforward extensions of the estimates of the lower
norms. It seems to be a general feature of the problem that, under the
natural additional regularity assumption on $(u_a, A_a)$ and/or on
$(u_+, A_+, \dot{A}_+)$, \par

(i) the required time decay remains the same, and in particular the
value of $\lambda$ needed for Step 1 remains $\lambda = 3/8$, \par

(ii) no additional smallness condition appears on the higher norms of
$u_a$ and/or $u_+$ beyond the smallness condition of $c_4$ that appears
at the lower level.\par

Two theories with higher level of regularity are of special interest~:
\par

(1) The theory at level $(k, \ell) = (1,1)$ since that is the level of
the energy for the (WS)$_3$ system. The appropriate function space is
\bea
\label{2.59e}
&&X_1(I) = \Big \{ (v, B):v\in {\cal C}(I, H^1), \ \nabla B, \partial_t B \in {\cal C}(I, L^2),\nn\\
&&\parallel (v, B);X_1(I)\parallel \ \equiv \ \mathrel{\mathop {\rm Sup}_{t \in I }}\ h(t)^{-1} \Big ( \parallel v(t);H^1\parallel \ + \ \parallel v; L^{8/3}(J,W_4^1)\parallel\nn \\
&& + \ \parallel B; L^4(J, L^4)\parallel\ + \ \parallel \nabla B(t) \parallel_2 \ + \ \parallel \partial_t B(t) \parallel_2 \Big ) < \infty \Big \} 
\eea

\noi where $J = [t, \infty ) \cap I$. It differs from the previous
space $X(I)$ by the inclusion of the $L^2$ norm and of a Strichartz
norm of $\nabla v$, and of the energy norm of $B$. Actually in that
definition, the $L^2$ norm of $\partial_tB$ is optional since it is
never used to perform the contraction estimates. It can be recovered at
the end as a by-product. Furthermore it turns out that the energy norm
of $B$ has a better time decay than that expressed by the definition of
$X_1(\cdot )$. The additional estimates needed for that theory as
compared with the lower level one are obtained from $L^2$ and
Strichartz estimates on the gradient of the Schr\"odinger equation,
namely on (\ref{2.4e}) with $\partial = \nabla$, and from the energy
estimate (\ref{2.7e}) for $B$. We refer to \cite{14r} Proposition 1.1,
part (2), Proposition 2.3 and Proposition 3.1, part (2) for a
description of that theory with the simplest asymptotic $(u_a, A_a)$.
\par

(2) The theory at level $(k, \ell ) = (2,1)$ since $k = 2$ is the
lowest level for $v$ where we can treat the (MS)$_3$ system. That theory
is therefore a simplified model for the (MS)$_3$ theory, which has a
strong similarity with it. The appropriate function space is now 
\bea
\label{2.60e}
&&X_2(I) = \Big \{ (v, B):v\in {\cal C}(I, H^2) \cap {\cal C}^1(I, L^2), \ \nabla B, \partial_t B \in {\cal C}(I, L^2),\nn\\
&&\parallel (v, B);X_2(I)\parallel\ \equiv \ \mathrel{\mathop {\rm Sup}_{t \in I }}\ h(t)^{-1} \Big ( \parallel v(t);H^2\parallel \ + \ \parallel \partial_t v(t) \parallel_2 \nn\\									 
&&+ \ \parallel v; L^{8/3}(J,W_4^2)\parallel \  + \ \parallel \partial_t v;  L^{8/3}(J,L^4)\parallel \ + \ \parallel B; L^4(J, L^4)\parallel\nn \\ 
&&+ \ \parallel \nabla B(t) \parallel_2 \ + \ \parallel \partial_t B(t) \parallel_2 \Big ) < \infty \Big \} 
\eea
 
\noi where $J = [t, \infty ) \cap I$. It differs from the previous
$X(\cdot )$ by the inclusion of the $L^2$ norm and of a Strichartz norm
of $\partial_tv$ and of $\Delta v$, and of the energy norm of $B$. Since
the Schr\"odinger equation directly relates $\Delta v$ to $\partial_t
v$, the norms of $\Delta v$, namely of level $k=2$, are estimated in
terms of those of $\partial_tv$. As a result the basic estimates of $v$
can be performed on the time derivative of the Schr\"odinger equation,
namely on (\ref{2.4e}) with $\partial = \partial_t$. It is then sufficient
to apply only one (time) derivative on $B$, so that $B$ can still be
taken at the level $\ell = 1$ only, as in the case of the $H^1$ theory.
Furthermore the $L^2$ norm of $\nabla B$ and the $L^2$ norm and
Strichartz norm of $\Delta v$ are optional in the definition of
$X_2(\cdot)$, since they are not used in the contraction proof to
estimate the other norms that occur in the definition. They can be
recovered at the end as a by-product. Finally, as in the case of the
$H^1$ theory, it turns out that the energy norm of $B$ has a better
time decay than that expressed by the definition of $X_2(\cdot )$.\par

In the remaining part of this section, we describe the results
available at the level $(k, \ell ) = (2,1)$, as a simplified model of
the (MS)$_3$ theory to be presented in the next section. We refer to
\cite{14r}, Proposition 1.1, part (3), Proposition 2.4 and Proposition
3.1, part (3) for a description of that theory with the simplest
asymptotic $(u_a, A_a)$. The result with improved asymptotics is
adapted from \cite{32r}.\par

The basic result concerning Step 1 is the following. \\

\noi {\bf Proposition 2.4.} {\it Let $h$ be defined as previously with
$\lambda = 3/8$ and let $X_2(\cdot )$ be defined by (\ref{2.60e}). Let
$(u_a, A_a)$ satisfy the estimates (\ref{2.13e})-(\ref{2.16e}) of
Proposition 2.1 and in addition
\beq
\label{2.61e}
\parallel u_a(t) \parallel_{\infty} \ \leq  c\ t^{-3/2} \qquad , \quad \parallel \partial_t u_a (t) \parallel_4\ \leq c\ t^{-3/4} \ ,
\eeq
\beq
\label{2.62e}
\parallel \partial_t A_a (t)\parallel_{\infty} \ \leq a \ t^{-1} \ ,
\eeq
\beq
\label{2.63e}
\parallel \partial_t R_1;L^1([t, \infty ), L^2 ) \parallel\ \leq r_1\ h(t) \ ,
\eeq
\beq
\label{2.64e}
\parallel R_1;L^{8/3} ([t, \infty ), L^4) \parallel\ \leq r_1\ h(t) \,
\eeq
\beq
\label{2.65e}
\parallel R_2;L^1 ([t, \infty ), L^2) \parallel\ \leq r_2\ t^{-1/2}\ h(t) 
\eeq

\noi for some constants $c_4$, $c$, $a$, $r_1$ and $r_2$ with $c_4$
sufficiently small and for all $t \geq 1$. Then there exists $T$, $1
\leq T < \infty$ and there exists a unique solution $(v, B)$ of the
system (\ref{1.6e}) in the space $X_2([T, \infty ))$. Furthermore $B$ satisfies the estimate
\beq
\label{2.66e}
\parallel \nabla B(t) \parallel_2\ \vee\ \parallel \partial_t B(t)\parallel_2\ \leq C \left ( t^{-1/2} + t^{1/4} h(t)\right ) h(t)
 \eeq

\noi for some constant $C$ and for all $t \geq T$.}\\

The proof follows closely that of Proposition 2.1. One starts from $(v,
B) \in X_2([T, \infty ))$ for some $T$, $1 \leq T < \infty$, and
therefore satisfying (\ref{2.17e}) and in addition
\beq
\label{2.67e}
\left \{ \begin{array}{l}
\parallel \partial_t v(t) \parallel_2 \ \leq N_3\ h(t)\\
\\
\parallel \partial_t v ; L^{8/3}([t, \infty ), L^4)\parallel \ \leq \  N_4\ h(t)\\ 
\\
\parallel \nabla B(t)\parallel_2 \ \vee \ \parallel \partial_t B(t)\parallel_2\ \leq \  N_5\ h(t) \\
\\
\parallel \Delta v(t) \parallel_2\ \leq N_6\ h(t) \\
\\
\parallel \Delta v; L^{8/3}([t, \infty ), L^4) \parallel\ \leq N_7\ h(t) \\
\end{array}\right . \eeq

\noi for some constants $N_i$, $0 \leq i \leq 7$ and for all $t \geq
T$. For each such $(v,B)$ one constructs a solution $(v', B')$ of the system (\ref{1.8e}) in
$X_2(I)$, first for $T = [T, t_0]$ and then for $I = [T, \infty )$, and
one shows finally that the map $\phi : (v, B) \to (v', B')$ thereby defined is a
contraction on the subset ${\cal R}_2$ of $X_2([T, \infty ))$ defined
by (\ref{2.17e}) (\ref{2.67e}) for suitably chosen $N_i$ and for
sufficiently large $T$. The crux of the proof is to estimate the seminorms (\ref{2.18e}) of $(v',B')$ and the seminorms $N'_i$, $3 \leq i
\leq 7$ defined in an obvious way in correspondence with (\ref{2.67e}).
The fact that no additional smallness assumption is needed to complete
the proof comes from the fact that the linear part of the system of
conditions extending (\ref{2.21e}) and including the additional
variables $N_i$, $3 \leq i \leq 7$, is essentially triangular with
respect to the latter variables.\par

We next state the final result that can be obtained with the simple asymptotics considered above.\\

\noi {\bf Proposition 2.5.} {\it Let $h(t) = t^{-1/2}$ and let
$X_2(\cdot )$ be defined by (\ref{2.60e}). Let $(u_a, A_a)$ be defined
by (\ref{2.29e}) (\ref{2.35e}) (\ref{2.23e}) (\ref{2.24e})
(\ref{2.30e}) (\ref{2.31e}). Let $u_+ \in H^{1,3} \cap H^{2,2}$ with
$c_4 = \ \parallel \widehat{u}_+ \parallel_4$ sufficiently small. Let $A_+$,
$\omega^{-1}\dot{A}_+\in H^1$ and $\nabla^2A_+$, $\nabla \dot{A}_+ \in
W_1^1$. Then there exists $T$, $1 \leq T < \infty$, and there exists a
unique solution $(u, A)$ of the (WS)$_3$ system (\ref{2.1e}) such that
$(v, B) \equiv (u- u_a, A- A_a) \in X_2 ([T, \infty ))$. Furthermore
$B$ satisfies the estimate
\beq
\label{2.68e}
\parallel \nabla B(t) \parallel_2\ \vee \ \parallel \partial_t B(t) \parallel_2\ \leq C\ t^{-3/4}
\eeq

\noi for some constant $C$ and for all $t \geq T$.}\\

We finally state the final result that can be obtained by using the improved asymptotics of \cite{32r}. \\

\noi {\bf Proposition 2.6.} {\it Let $h(t) = t^{-1}(2 + \ell n\ t)^2$ and let
$X_2(\cdot )$ be defined by (\ref{2.60e}). Let $(u_a, A_a)$ be defined
by (\ref{2.46e}) (\ref{2.47e}) (\ref{2.49e}) and (\ref{2.35e})
(\ref{2.23e}) (\ref{2.24e}) (\ref{2.30e}) (\ref{2.31e}). Let $u_+ \in H^{1,3} \cap H^{2,2}$ with
$xu_+ \in W_1^2$, $x^2 u_+ \in W_1^1$ and with $c_4 = \ \parallel \widehat{u}_+ \parallel_4$ sufficiently small. Let $(A_+, \dot{A}_+)$ satisfy
\bea
\label{2.69e}
&&A_+, \omega^{-1} \dot{A}_+ \in \dot{H}^{-2} \cap H^1 \quad , \quad \nabla^2 A_+, \nabla \dot{A}_+ \in W_1^1\ ,\nn \\
&&x\cdot \nabla A_+, \omega^{-1} x \cdot \nabla \dot{A}_+ \in \dot{H}^{-2} \cap \dot{H}^{-1}\ .
\eea

\noi Then :\par

(1)  There exists $T$, $1 \leq T < \infty$, and there exists a
unique solution $(u, A)$ of the (WS)$_3$ system (\ref{2.1e}) such that
$(v, B) \equiv (u- u_a, A- A_a) \in X_2 ([T, \infty ))$. Furthermore
$B$ satisfies the estimate
\beq
\label{2.70e}
\parallel \nabla B(t) \parallel_2\ \vee \ \parallel \partial_t B(t) \parallel_2\ \leq C\ t^{-3/2} (2 + \ell n\ t)^2
\eeq

\noi for some constant $C$ and for all $t \geq T$.\par

(2) There exists $T$, $1 \leq T < \infty$, and there exists a
unique solution $(u, A)$ of the (WS)$_3$ system (\ref{2.1e}) such that
$(u- u_1, A- A_a) \in X_2 ([T, \infty ))$. One can take the same $T$ and the solution $(u, A)$ is the same as in Part (1).}\\

As in Proposition 2.3, Part (2) follows immediately from Part (1) and
from the fact that $(fu_1,0) \in X_2([T, \infty ))$ for any $T \geq 1$
under the assumptions made on $(u_+, A_+, \dot{A}_+)$.

\mysection{The Maxwell-Schr\"odinger system (MS)$_3$}
\hspace*{\parindent}
In this section we review the main results available on the local
Cauchy problem at infinity for the (MS)$_3$ system. That system can be
written as 
\beq
\label{3.1e} \left \{ \begin{array}{l}  i\partial_t u = - (1/2)
\Delta_A u + A_e u \\ \\ \sq A_e  - \partial_t \left ( \partial_t A_e
+ \nabla \cdot A\right ) = |u|^2  \\ \\ \sq A + \nabla \left (
\partial_t A_e + \nabla \cdot A\right ) = {\rm Im} \ \overline{u}
\nabla_A u \ .  \end{array} \right . \eeq

\noi Here $(A,A_e)$ is an ${I\hskip-1truemm R}^{3+1}$ valued function
defined in space time ${I\hskip-1truemm R}^{3+1}$, $\nabla_A = \nabla -
i A$ and $\Delta_A = \nabla_A^2$ are the covariant gradient and
Laplacian respectively. An important property of that system is its
gauge invariance, namely the invariance under the transformation
\beq
\label{3.2e}
(u, A, A_e) \to \left ( u \exp (§-i \varphi ), A - \nabla \varphi, A_e + \partial_t \varphi \right ) \ ,
\eeq

\noi where $\varphi$ is an arbitrary real function defined in
${I\hskip-1truemm R}^{3+1}$. As a consequence of that invariance, the
system (\ref{3.1e}) is underdetermined as an evolution system and has
to be supplemented by an additional equation, called a gauge condition.
Here we shall use the Coulomb gauge condition, namely $\nabla \cdot A =
0$. Under that condition, one can replace the system (\ref{3.1e}) by a
formally equivalent one in the following standard way. The second
equation of (\ref{3.1e}) can be solved for $A_e$ by 
\beq
\label{3.3e}
A_e = - \Delta^{-1} |u|^2 = ( 4\pi |x|)^{-1} \ * \ |u|^2 \equiv g(|u|^2) \ .
\eeq

\noi Substituting (\ref{3.3e}) into the first and last equation of (\ref{3.1e}) yields the new system
\beq
\label{3.4e}
\left \{ \begin{array}{l}  i\partial_t u = - (1/2)
\Delta_A u + g(|u|^2) u \\ \\ \sq A  = P\ {\rm Im} \ \overline{u} \nabla_A u  \qquad , \quad \nabla \cdot A = 0\end{array} \right . 
\eeq

\noi where $P = \1 - \nabla \Delta^{-1} \nabla$ is the projector on
divergence free vector fields so that the gauge condition $\nabla\cdot
A = 0$ is preserved by the evolution. The system (\ref{3.4e}) coincides
with that given in the introduction.\par

The Coulomb gauge has several advantages. In the expansion of the
covariant Laplacian, the term $iA\cdot \nabla$ is a transport
term by the vector field $A$. Using the Coulomb gauge consists in
decomposing $A$ into a divergence free part and a gradient and taking
advantage of the gauge invariance of the system to eliminate the
gradient part. It makes the transport term isometric in $L^r$ for any
$r$ and eliminates all derivatives of $A$ in the system. The Coulomb
gauge is also strongly preferred by atomic physicists since it provides
a clean separation of the electromagnetic interaction into an
instantaneous electrostatic part and a propagating magnetic part. \par

Nevertheless other gauges can be considered for the system
(\ref{3.1e}). The Lorentz gauge is defined by $\partial_t A_e + \nabla
\cdot A = 0$. In that gauge, the (MS)$_3$ system becomes
\beq
\label{3.5e} \left \{ \begin{array}{l}  i\partial_t u = - (1/2)
\Delta_A u + A_e u \\ \\ \sq A_e  = |u|^2  \\ \\ \sq A  = {\rm Im} \ \overline{u}
\nabla_A u \\ \\
\partial_t A_e + \nabla \cdot A = 0 \ .  \end{array} \right . \eeq

\noi In particular the electric potential $A_e$ is now coupled to $u$
by a WS type interaction, with however a change of sign. As in the case
of the Coulomb gauge, the gauge condition is preserved by the
evolution.\par

The Lorentz gauge condition is Lorentz invariant, which is of little
advantage for the MS system since the Schr\"odinger equation is not.
That invariance becomes important for Lorentz invariant systems such as
the Maxwell-Dirac or Maxwell-Klein-Gordon systems, for instance in
order to ensure that Lorentz invariance is preserved term by term in
formal expansions such as perturbation expansions. This property is
extensively used for instance in Quantum Electrodynamics. \par

In this section, we shall use the Coulomb gauge and make only a brief
comment at the end on the Lorentz gauge.\par

In contrast with the case of the other equations and systems listed in
the introduction, the Cauchy problem for the (MS)$_3$ system is not in a
satisfactory shape. It is known that this problem is locally well posed
at the level $(k, \ell )$ for sufficiently large $(k, \ell )$, namely
with sufficiently high regularity \cite{21r} \cite{22r}, the best
result available so far being $(k, \ell ) = (5/3, 4/3)$ \cite{22r}. On
the other hand, the (MS)$_3$ system has global solutions at the level
of the energy, namely $(k, \ell ) = (1, 1)$, obtained by compactness
and therefore without uniqueness \cite{17r}. However the (MS)$_3$
system is so far not known to be globally well posed in any function
space. \par

We now turn to the local Cauchy problem at infinity for the (MS)$_3$
system in the Coulomb gauge (\ref{3.4e}). The exposition follows
\cite{15r} which is strongly inspired by \cite{36r}, and includes also
some results from \cite{33r}. Note that \cite{33r} \cite{36r} also
consider the (MS)$_3$ system in the Lorentz gauge (\ref{3.5e}).\par

The theory follows the same pattern as that for the (WS)$_3$ system and
strongly resembles the $H^2$ theory for the latter. One looks for $(u,
A)$ in the form $(u, A) = (u_a + v, A_a + B)$ where $(u_a, A_a)$ is the
asymptotic form, which has to satisfy the gauge condition $\nabla\cdot
A_a = 0$. The auxiliary system satisfied by the new functions $(v, B)$
is now
\beq \label{3.6e}
\left \{ \begin{array}{l}  i\partial_t v = - (1/2)
\Delta_A v + g(|u|^2)v + G_1 - R_1 \\ \\ \sq B =  G_2 - R_2 \qquad , \quad \nabla \cdot B = 0 \ ,  \end{array} \right . 
 \eeq

\noi where $G_1$ and $G_2$ are defined by
\beq \label{3.7e}
\left \{ \begin{array}{l}  G_1= iB \cdot \nabla_{A_a} u_a + (1/2)
B^2 u_a + g\left ( |v|^2 + 2\ {\rm Re}\ \overline{u}_a v \right ) u_a \\ \\ G_2 = P\ {\rm Im}\left ( \overline{v} \nabla_A v + 2 \overline{v}\nabla_A u_a\right )- P\ B |u_a|^2  \end{array} \right . 
 \eeq

\noi and the remainders are defined by
\beq \label{3.8e}
\left \{ \begin{array}{l} R_1 =  i\partial_t u_a + (1/2)
\Delta_{A_a} u_a - g \left ( |u_a|^2\right ) u_a \\ \\ R_2 = \sq A_a  - P \ {\rm Im} \ \overline{u}_a \nabla_{A_a} u_a \ . \end{array} \right . 
 \eeq

\noi As in the case of (WS)$_3$, we consider also the partly linearized system for functions $(v', B')$
\beq \label{3.9e}
\left \{ \begin{array}{l}  i\partial_t v' = - (1/2)
\Delta_A v' + g(|u|^2)v' + G_1 - R_1 \\ \\ \sq B' =  G_2 - R_2 \ , \qquad \nabla\cdot B' = 0   \ .\end{array} \right . 
 \eeq
 
\noi The first step of the method consists again in solving the system
(\ref{3.6e}) for $(v, B)$, with $(v, B)$ tending to zero at infinity in
time in suitable norms, under assumptions on $(u_a, A_a)$ of a general
nature, the most important of which being decay assumptions on the
remainders $R_1$ and $R_2$.

When choosing the appropriate function space $X(\cdot )$ to perform that step, one has to take into account the following difference with the (WS)$_3$ case. If $v$ satisfies the Schr\"odinger equation 
\beq
\label{3.10e}
i \partial_t v = - (1/2) \Delta_A v + Vv + f 
\eeq

\noi with real $V$, then again
$$\partial_t \parallel v\parallel_2^2 \ = 2\ {\rm Im}\ <v, f>$$

\noi and $v$ is estimated in $L^2$ by (\ref{2.3e}) independently of $A$
and $V$. However in order to use the Strichartz inequality
(\ref{2.9e}), one has to expand the covariant Laplacian, so that 

\beq
\label{3.11e}
i\partial_t v = - (1/2) \Delta v + i A \cdot \nabla v + \ \hbox{other terms}
\eeq

\noi and the term $A \cdot \nabla v$ has to be included in $f$, thereby
requiring the control of one more derivative than appears at the level
where the Strichartz inequality is considered. As a consequence, when
using $L^2$ norms at the level $k$, one can use Strichartz norms only
at the level $k - 1$. \par

Taking that point into account, one can define the relevant space $X(\cdot )$ by 
$$X(I) = \Big \{ (v, B):v\in {\cal C}(I, H^2) \cap {\cal C}^1(I,L^2),$$
$$\parallel (v, B);X(I)\parallel\ \equiv \ \mathrel{\mathop {\rm Sup}_{t \in I }}\ h(t)^{-1} \Big ( \parallel v(t);H^2\parallel \ + \ \parallel \partial_t v(t) \parallel_2 \ + \ \parallel v; L^{8/3}(J,W_4^1)\parallel$$
\beq
\label{3.12e}
+ \ \parallel B;L^4(J,W_4^1)\parallel \ + \ \parallel \partial_t B;L^4 (J, L^4)\parallel \Big ) < \infty \Big \} \eeq

\noi where $J = [t, \infty) \cap I$ and $h$ is defined as in Section 2.
That space has level $(k, \ell ) = (2, 3/2)$. More precisely it
includes $L^2$ norms of $v$ at the level $k=2$, Strichartz norms of $v$
at the level $k=1$, in keeping with the preceding remark, and
Strichartz norms of $B$ at the level $\ell = 3/2$ and $\ell = 1/2$
(compare with Lemma 2.2, especially (\ref{2.10e}) (\ref{2.11e})). One
could also include the energy norm of $B$, which is at the level $\ell
= 1$, but that norm is not needed for the proof of the estimates and it
comes out at the end as a by-product thereof (see Propositions 3.1 and
3.2 below). The space (\ref{3.12e}) is larger than that used in
\cite{33r} \cite{36r} which has level $(k, \ell ) = (3, 2)$. \par

We can now state the main result concerning Step 1, which is Proposition 2.2 of \cite{15r}.\\

\noi {\bf Proposition 3.1.} {\it Let $h$ be defined as in Section 2 with
$\lambda = 3/8$ and let $X(\cdot )$ be defined by (\ref{3.12e}). Let
$u_a$, $A_a$, $R_1$ and $R_2$ be sufficiently regular (for the following estimates to make sense) and satisfy the estimates
\beq
\label{3.13e}
\parallel \partial_t^j \nabla^k u_a(t) \parallel_r \ \leq \ c\ t^{-\delta (r)} \quad \hbox{\it for $2 \leq r \leq \infty$}
\eeq

\noi where $\delta (r) = 3/2 - 3/r$, and in particular
\beq
\label{3.14e}
\parallel u_a\parallel_3\ \leq c_3 \ t^{-1/2}\quad , \quad \parallel \nabla u_a\parallel_4\ \leq c_4 \ t^{-3/4}\ ,
\eeq
\beq
\label{3.15e}
\parallel \nabla^2 u_a(t) \parallel_{4} \ \vee \ \parallel \partial_t \nabla u_a (t) \parallel_{4}\ \leq c\ t^{-3/4}\ ,
\eeq
\beq
\label{3.16e}
\parallel \partial_t^j \nabla^k A_a(t) \parallel_{\infty} \ \leq a\ t^{-1}\ ,
\eeq
\beq
\label{3.17e}
\parallel \partial_t^j \nabla^k R_1;L^1([t , \infty ),L^2) \parallel \ \leq r_1\ h(t)\ ,
\eeq
\beq
\label{3.18e}
\parallel R_2;L^{4/3}([t, \infty ), W_{4/3}^1) \parallel\  \leq r_2\ h(t)\ ,
\eeq

\noi for $0 \leq j+k \leq 1$, for some constants $c$, $c_3$, $c_4$, $a$, $r_1$ and $r_2$ with $c_3$ and $c_4$
sufficiently small and for all $t \geq 1$. Then there exists $T$, $1 \leq T < \infty$
and there exists a unique solution $(v, B)$ of the
system (\ref{3.6e}) in $X([T, \infty ))$. If in addition 
\beq
\label{3.19e}
\parallel R_2;L^{1}([t, \infty ), L^{2}) \parallel\ \leq r_2\ t^{-1/2}\ h(t)\ ,
\eeq

\noi for all $t \geq T$, then $\nabla B$, $\partial_t B \in {\cal
C}([T, \infty ), L^2)$ and $B$ satisfies the estimate 
\beq
\label{3.20e}
\parallel \nabla B(t)\parallel_2\ \vee \ \parallel \partial_t B(t)\parallel_2\ \leq C \left (  t^{-1/2}  + t^{1/4}\ h(t) \right )  h(t)
\eeq

\noi for some constant $C$ and for all $t \geq T$.}\\

The proof follows closely those of Propositions 2.1 and 2.4. One starts
from $(v,B)\in X([T, \infty ))$ for some $T$, $1 \leq T < \infty$, so
that $(v, B)$ satisfies
\beq
\label{3.21e}
\left \{ \begin{array}{l} \parallel v(t) \parallel_2 \ \leq N_0\ h(t) \\
\\
\parallel v ; L^{4}(J, L^3) \parallel \ \vee \ \parallel v ; L^{8/3}(J, L^4) \parallel\ \leq N_1\ h(t)\\
\\
\parallel B ; L^{4}(J, L^4) \parallel \ \leq N_2\ h(t)\\
\\
\parallel \partial_t v(t) \parallel_2 \ \leq N_3\ h(t)\\
\\
\parallel \nabla v ; L^{4}(J, L^3) \parallel \ \vee \ \parallel \nabla v ; L^{8/3}(J, L^4) \parallel\ \leq N_4\ h(t)\\
\\
\parallel \Delta v(t) \parallel_2 \ \leq N_5\ h(t)\\
\\
\parallel \nabla B ; L^{4}(J, L^4) \parallel \ \vee \ \parallel \partial_t B  ; L^{4}(J, L^4) \parallel\ \leq N_6\ h(t)
\end{array} \right .
\eeq

\noi for some constants $N_i$, $0 \leq i \leq 6$ and for all $t \geq
T$, with $J = [t, \infty )$. For each such $(v, B)$ one constructs a solution $(v', B')$
of the system (\ref{3.9e}) in $X([T, \infty ))$ and one shows that the
map $\phi : (v, B) \to (v', B')$ thereby defined is a contraction on
the subset ${\cal R}$ of $X([T, \infty ))$ defined by (\ref{3.21e}) for
suitably chosen $N_i$ and for sufficiently large $T$. The smallness
conditions on $c_3$ and $c_4$ appear when trying to solve the extension
to the present case of the system (\ref{2.21e}). The coefficient $c_3$
comes from the Hartree interaction and appears in a diagonal term in the
equation for $N_0$. The coefficient $c_4$ comes from the magnetic
interaction and appears in cross couplings of $N_0$ and $N_2$ as in
(\ref{2.21e}). The rest of the system is triangular in its linear part
and does not generate any additional smallness condition. Note also
that the decay exponent $\lambda = 3/8$ is the same as for the (WS)$_3$
system, and that the energy norm of $B$ has a stronger decay than that occurring
in the definition of $X(\cdot )$. \par

We now turn to Step 2, namely to the construction of $(u_a, A_a)$
satisfying the assumptions of Proposition 3.1. Again solutions of the
underlying free equations are immediately seen to be inadequate, but
the weak time decay allowed by Proposition 3.1 allows for the simplest
modification available in the literature. We want to ensure sufficient
decay properties of the remainders. We consider first $R_2$. Using the
operator
\beq
\label{3.22e}
J = x + it\nabla
\eeq

\noi we rewrite $R_2$ as
\beq
\label{3.23e}
R_2 = \sq A_a + P \left ( t^{-1} \ {\rm Re}\ \overline{u}_a \ J u_a + A_a |u_a|^2 - t^{-1}x |u_a|^2 \right ) \ .
\eeq

\noi We now choose $A_a$ according to (\ref{2.23e}) (\ref{2.24e}), now however with 
\beq
\label{3.24e}
A_1 = \int_t^{\infty} dt' (\omega t')^{-1} \sin (\omega (t'-t))P\ x |u_a (t')|^2 \ ,
\eeq

\noi thereby ensuring that $\sq A_0 = 0$, $\sq A_a = \sq A_1 = Pt^{-1} x |u_a|^2$ and therefore
\beq
\label{3.25e}
R_2 = P \left ( t^{-1} \ {\rm Re}\ \overline{u}_a\ J u_a + A_a |u_a|^2 \right ) \ .
\eeq

\noi With that choice, under general assumptions on $(u_a, A_a)$ of the
same type as in Proposition 3.1 but not making use of their special
form, one can prove that $R_2$ satisfies the assumptions needed for
that proposition with $h(t) = t^{-1} (2 + \ell n\ t)$.\\

\noi {\bf Lemma 3.1.} {\it Let $(u_a, A_a)$ satisfy the estimates
\beq
\label{3.26e}
\parallel u_a(t);W_r^1\parallel \ \leq c\ t^{-\delta (r)} \qquad \hbox{\it for $2 \leq r \leq \infty$} \ ,
\eeq
\noi where $\delta (r) = 3/2 - 3/r$,
\beq
\label{3.27e}
\parallel J u_a(t)\parallel_2 \ \leq c_1\ (2 + \ell n\ t) \ ,
\eeq
\beq
\label{3.28e}
\parallel A_a(t);W_{\infty}^1\parallel\ \leq a\ t^{-1} 
\eeq

\noi for all $t \geq 1$. Then $R_2$ satisfies the estimates
\beq
\label{3.29e}
\parallel R_2;L^{4/3}([t, \infty ),W_{4/3}^1)\parallel \  \leq r_2\ t^{-1}(2 + \ell n \ t) \ , 
\eeq
\beq
\label{3.30e}
\parallel R_2;L^1([t, \infty ),L^2)\parallel \ \leq r_2\ t^{-3/2} (2 + \ell n \ t) 
\eeq

\noi for some constant $r_2$ and for all $t \geq 1$.}\\

We now turn to $R_1$. We first skim $R_1$ of some harmless terms.
Expanding the covariant Laplacian and using again $J$, we rewrite
$R_1$ as
\beq
\label{3.31e}
R_1 = R_{1,1} + R_{1,2}
\eeq

\noi where
\beq
\label{3.32e}
R_{1,1} = i \partial_t u_a + (1/2) \Delta u_a + t^{-1} (x \cdot A_1) u_a - g(|u_a|^2)u_a \ ,
\eeq
\beq
\label{3.33e}
R_{1,2} = t^{-1}(x \cdot A_0) u_a - t^{-1} A_a \cdot Ju_a - (1/2) A_a^2 u_a \ .
\eeq

\noi In the same way as for $R_2$, one can show that $R_{1,2}$ satisfies
the assumptions needed for Proposition 3.1 with $h(t) = t^{-1} (2 + \ell n\ t)$
under general assumptions on $(u_a, A_a)$
not making use of their special form.\\

\noi {\bf Lemma 3.2.} {\it Let $u_a$, $A_a$ and $A_0$ satisfy the estimates
\beq
\label{3.34e}
\parallel \partial_t^j \nabla^k u_a \parallel_2\ \leq c \ ,
\eeq
\beq
\label{3.35e}
\parallel \partial_t^j \nabla^k J u_a \parallel_2\ \leq c_1(2 + \ell n\ t) \ ,
\eeq
$$\parallel \partial_t^j \nabla^k A_a \parallel_{\infty}\ \leq a\ t^{-1} \ , \eqno(3.16)\equiv(3.36)$$
$$\parallel \partial_t^j \nabla^k (x \cdot A_0) \parallel_{\infty}\ \leq a_0\ t^{-1} \ , \eqno(3.37)$$

\noi for $0 \leq j + k \leq 1$ and for all $t \geq 1$. Then $R_{1,2}$ satisfies the estimates 
$$\parallel \partial_t^j \nabla^k R_{1,2} \parallel_{2}\ \leq r_{1,2} \ t^{-2} (2 + \ell n \ t)\ , \eqno(3.38)$$

\noi for $0 \leq j + k \leq 1$, for some constant $r_{1,2}$ and for all $t \geq 1$.}\\

We next choose $u_a$ according to (\ref{2.29e}). Substituting
(\ref{2.29e}) into (\ref{3.24e}) yields again (\ref{2.30e}), now however
with 
$$\widetilde{A}_1 =  \int_1^{\infty} d\nu \ \nu^{-3} \omega^{-1} \sin (\omega (\nu - 1)) D_0 (\nu ) P\ x|\widehat{u}_+|^2 \ .
\eqno(3.39)$$

Again $\widetilde{A}_1$ is constant in time. Substituting (\ref{2.29e})
(\ref{2.30e}) into the definition (\ref{3.32e}) of $R_{1,1}$ and using
again the commutation relation (\ref{2.32e}), we obtain
$$R_{1,1} = MD \left ( i \partial_t + (2t^2)^{-1} \Delta + t^{-1} x \cdot \widetilde{A}_1 - t^{-1} g(|\widehat{u}_+|^2)\right ) \exp (-i \varphi ) \widehat{u}_+ \ . \eqno(3.40)$$

\noi We finally choose $\varphi$ so as to cancel the long range terms
$x\cdot \widetilde{A}_1$ and $g(|\widehat{u}_+|^2)$ in (3.40). We take
$$\varphi = (\ell n\ t) \left ( g(|\widehat{u}_+|^2) - x\cdot \widetilde{A}_1 \right )\ ,
\eqno(3.41)$$

\noi thereby obtaining
$$R_{1,1} = (2t^2)^{-1} \ MD\ \Delta \exp (- i \varphi ) \widehat{u}_+ \eqno(3.42)$$

\noi which is very similar to (\ref{2.36e}) except for the different choice of $\varphi$.\par

Using mainly Sobolev inequalities, one can derive the estimates of
$(u_a, A_1)$ needed for Lemmas 3.1 and 3.2, and the remaining estimates
of $u_a$ and of $R_{1,1}$ needed for Proposition 3.1. The following
lemma is slightly stronger than needed (see Proposition 3.1 in
\cite{15r}). \\

\noi {\bf Lemma 3.3.} {\it Let $u_a$ be defined by (\ref{2.29e})
with $\varphi$ defined by (3.41)
(\ref{3.3e}) (3.39) and let $R_{1,1}$ be given by (3.42). Let
$u_+ \in H^{1,3} \cap H^{3,1}$. Then the following estimates hold for
some constants $c$, $c_1$, $a_1$ and $r_{1,1}$, for $0 \leq j+k \leq 1$ and
for all $t\geq 1$
$$\parallel \partial_t^j \nabla^k u_a(t)\parallel_r \ \leq c\ t^{-\delta (r)} \qquad \hbox{\it for $2 \leq r \leq \infty$}\ . \eqno(3.13)\equiv(3.43)$$

\noi In particular}
$$\parallel u_a(t) \parallel_3\ \leq \ \parallel \widehat{u}_+ \parallel_3 \ t^{-1/2} \ , \eqno(3.44)$$
$$\parallel \nabla u_a(t) \parallel_4\ \leq \left (  \parallel x\widehat{u}_+ \parallel_4\ + O(t^{-1} \ell n\ t)\right ) t^{-3/4} \ .   \eqno(3.45)$$
$$\parallel \partial_t^j \nabla^{k+1} u_a(t) \parallel_r\ \leq \ c\ t^{-\delta (r)}\qquad \hbox{\it for $2 \leq r \leq 6$} \ , \eqno(3.46)$$
$$\parallel \partial_t^j \nabla^{k} J u_a(t) \parallel_r\ \leq \ c_1(2 + \ell n\ t)\ t^{-\delta (r)}\qquad \hbox{\it for $2 \leq r \leq 6$} \ , \eqno(3.47)$$
$$\parallel \partial_t^j \nabla^{k} A_1 \parallel_{\infty}\ \leq \ a_1\ t^{-1} \ ,  \eqno(3.48)$$
$$\parallel \partial_t^j \nabla^{k} R_{1,1}(t) \parallel_2\ \leq \ r_{1,1}\ t^{-2}(2 + \ell n\ t)^2 \ .  \eqno(3.49)$$

\vskip 5 truemm
Finally, the contribution of $A_0$ to the conditions of Proposition 3.1
and of Lemmas 3.1 and 3.2, namely to (\ref{3.16e}) and (3.37), can be
estimated by the use of Lemma 2.3 and of the remark that if $A_0$ is a
solution of the wave equation $\sq A_0 = 0$ satisfying the Coulomb gauge
condition $\nabla\cdot A_0 = 0$, then $x\cdot A_0$ is also a solution
of the wave equation, namely $\sq (x\cdot A_0) = 0$. \par

Combining the previous estimates, especially Lemmas 3.1-3 with
Proposition 3.1, we obtain the final result for the (MS)$_3$ system in
the Coulomb gauge in the following form (see Proposition 1.1 in
\cite{15r}).\\

\noi{\bf Proposition 3.2.} {\it Let $h(t) = t^{-1}(2 + \ell n \ t
)^2$ and let $X(\cdot )$ be defined by (\ref{3.12e}). Let $u_a$ be
defined by (\ref{2.29e}) with $\varphi$ defined
by (3.41) (\ref{3.3e}) (3.39). Let $A_a = A_0 + A_1$
with $A_0$ defined by (\ref{2.24e}) and $A_1$ by (\ref{2.30e})
(3.39). Let $u_+ \in H^{3,1}\cap H^{1,3}$ with $\parallel x\widehat{u}_+
\parallel_4$ and $\parallel \widehat{u}_+ \parallel_3$ sufficiently small. Let
$\nabla^2 A_+$, $\nabla \dot{A}_+$, $\nabla^2(x\cdot A_+)$ and $\nabla
(x \cdot \dot{A}_+) \in W_1^1$ with $A_+$, $x \cdot A_+ \in L^3$ and
$\dot{A}_+$, $x \cdot \dot{A}_+ \in L^{3/2}$ and let $\nabla \cdot A_+
= \nabla \cdot \dot{A}_+ = 0$. \par

Then there exists $T$,
$1 \leq T < \infty$ and there exists a unique solution $(u,A)$ of the  (MS)$_3$
system (\ref{3.4e}) such that $(v, B) \equiv (u - u_a, A-A_a) \in X([T, \infty
))$. Furthermore $\nabla B$, $\partial_t B \in {\cal C}([T, \infty ), L^2)$ and $B$ satisfies the estimate
$$\parallel \nabla B(t) \parallel_2\ \vee \  \parallel \partial_t B(t) \parallel_2\ \leq C\ t^{-3/2} (2 + \ell n\ t)^2 \eqno(3.50)$$

\noi for some constant $C$ depending on $(u_+, A_+, \dot{A}_+)$ and for all $t \geq T$.}\\

\noi {\bf Remark 3.1} The only smallness conditions bear on $\parallel x \widehat{u}_+ \parallel_4$ and
on $\parallel \widehat{u}_+ \parallel_3$ and are required by the magnetic interaction and the Hartree interaction (\ref{3.3e}) respectively. 
In particular there is no smallness condition on $(A_+, \dot{A}_+)$.\\

\noi {\bf Remark 3.2.} The assumptions $A_+$, $x\cdot A_+ \in L^3$ and
$\dot{A}_+$, $x \cdot \dot{A}_+ \in L^{3/2}$ serve to exclude the
occurrence of constant terms in $A_+$, $x \cdot A_+$, $\dot{A}_+$, $x
\cdot \dot{A}_+$ and of terms linear in $x$ in $A_+$, $x \cdot A_+$,
but are otherwise implied by the $W_1^1$ assumptions on those
quantities through Sobolev inequalities.\\

\noi {\bf Remark 3.3.} The assumptions on $A_+$, $\dot{A}_+$ imply that
$\omega^{1/2}A_+$, $\omega^{-1/2} \dot{A}_+ \in H^1$ through Sobolev
inequalities. As a consequence the free wave solution $A_0$ defined by
(\ref{2.24e}) belongs to $L^4({I\hskip-1truemm R}, W_4^1)$ by Strichartz inequalities,
with $\partial_t A_0 \in L^4({I\hskip-1truemm R},L^4)$. In particular $A_0$
satisfies the local in time regularity of $B$ required in the
definition of the space $X(\cdot)$. Furthermore $\nabla A_+$, $\dot{A}_+ \in
L^2$ and therefore $\nabla A_0$,
$\partial_t A_0 \in ({\cal C} \cap L^{\infty})({I\hskip-1truemm R},L^2)$, namely $A_0$ is
a finite energy solution of the wave equation.\\

We conclude this section with some remarks on the (MS)$_3$ system in
the Lorentz gauge (\ref{3.5e}). There is every reason to expect similar
results in that case. Step 1 of the method should be implementable with
no smallness condition on $(A,A_e)$ and with the same decay exponent
$\lambda = 3/8$ as in Proposition 3.1. This should lead to a final
result with $h(t) = t^{-1/2}$ with simple asymptotics and with $h(t) =
t^{-1} (2 + \ell n\ t)^2$ with improved asymptotics, the correcting
term to $u_a$ being needed as in \cite{33r} to cancel the contribution
of $\nabla\cdot A_0$ coming from $\Delta_A$. The required level of
regularity cannot be expected to be lower than $(k, \ell, \ell_e) = (2,
2, 1)$ for $(v, B, B_e)$, with $\ell = 2$ for $B$ required by the
$\partial_t \nabla\cdot B$ term in the equation for $\partial_tv$ and
$\ell_e = 1$ for $B_e$ suggested by the (WS)$_3$ result of Proposition
2.4. Ascertaining whether that level is sufficient would require to go
through the detail of the estimates and we leave that as an open
question. In \cite{33r} \cite{36r} the Lorentz gauge case is treated
together with the Coulomb gauge case at the level $(k, \ell, \ell_e) =
(3,2,2)$.

\mysection{The Zakharov system (Z)$_{\bf n}$ for n = 2, 3}
\hspace*{\parindent} In this section we review the main results available on the local
Cauchy problem at infinity for the (Z)$_n$ system 
\beq
\label{4.1e}
\left \{ \begin{array}{l} i \partial_t u = - (1/2) \Delta u + A u\\ \\ \sq  A = \Delta  |u|^2\end{array} \right .
\eeq

\noi in space dimension $n=3$ and 2 (in that order). That system is
known to be globally well posed in the energy space (see \cite{2r}
\cite{6r} \cite{26r} and references to previous works therein quoted).
The exposition is based mostly on \cite{16r} and includes some previous
results from \cite{27r} \cite{34r}. \par

We first consider the case of dimension $n=3$. The (Z)$_3$ system is
short range, as will be clear below. We follow the sketch given in the
introduction and we first consider Step 1. For a given asymptotic
$(u_a, A_a)$, we look for $(u,A)$ in the form $(u, A) = (u_a + v,A_a+
B)$. The auxiliary system satisfied by the new functions $(v, B)$ is
now 
 \beq
\label{4.2e}
\left \{ \begin{array}{l} i \partial_t v = - (1/2) \Delta v + Av + Bu_a - R_1\\ \\ \sq  B = \Delta  \left ( |v|^2 + 2\ {\rm Re}\ \overline{u}_a v\right ) - R_2\end{array} \right .
\eeq

\noi where the remainders are defined by
\beq
\label{4.3e}
\left \{ \begin{array}{l} R_1 = i \partial_t u_a + (1/2) \Delta u_a - A_a u_a\\ \\ R_2 = \sq  A_a - \Delta  |u_a|^2\ .\end{array} \right .
\eeq

\noi As in the case of (WS)$_3$ we consider also the partly linearized system for functions $(v', B')$
 \beq
\label{4.4e}
\left \{ \begin{array}{l} i \partial_t v' = - (1/2) \Delta v' + Av' + Bu_a - R_1\\ \\ \sq  B' = \Delta  \left ( |v|^2 + 2\ {\rm Re}\ \overline{u}_a v\right ) - R_2\ . \end{array} \right .
\eeq

\noi Step 1 of the method consists again in solving the system
(\ref{4.2e}) for $(v, B)$ with $(v, B)$ tending to zero at infinity
under assumptions on $(u_a, A_a)$ of a general nature. In contrast with
the case of the (WS)$_3$ and (MS)$_3$ systems however, the Strichartz
inequalities for the wave equation do not seem to be useful for that
purpose, the reason being that the Laplacian in the RHS of the equation
for $B$ allows for an easy estimate of $B$ in $L^2$. The relevant space
$X(\cdot )$ can be defined as follows
\bea
\label{4.5e}
&&X(I) = \Big \{ (v, B):(v,B)\in {\cal C}(I, H^2\oplus H^1) \cap {\cal C}^1(I, L^2 \oplus L^2), \nn \\
&&\parallel (v, B);X(I)\parallel \ \equiv \ \mathrel{\mathop {\rm Sup}_{t \in I }}\ h(t)^{-1} 
\Big ( \parallel v(t);H^2\parallel \ + \ \parallel \partial_t v(t)\parallel_2\nn \\
&&+\ \parallel v; L^{8/3}(J,W_4^2)\parallel \ + \ \parallel \partial_t v; L^{8/3}(J, L^4)\parallel \nn\\
&& + \ \parallel B(t); H^1\parallel\ + \ \parallel \partial_t B(t) \parallel_2 \Big ) < \infty \Big \} 
\eea

\noi where $J = [t, \infty ) \cap I$ and $h$ is defined as in Section
2. That space has $(k, \ell ) = (2,1)$, namely 1 above the level of the
energy space in both variables. The definition includes $L^2$ norms and
Strichartz norms at the level $k = 2$ for $v$, and $L^2$ norms at the
levels $\ell = 0$ and $\ell = 1$ for $B$. We could also have included
the norm $\parallel \omega^{-1} \partial_t B\parallel_2$ which is part
of the Zakharov energy for $B$, but that norm is never used to perform the
estimates and it comes out at the end as a by-product thereof. We have
omitted it for simplicity.\par

We can now state the result concerning Step 1 (see Proposition 2.1 in
\cite{16r}). We recall that $\omega = (- \Delta)^{1/2}$.\\

\noi {\bf Proposition 4.1} {\it Let $h$ be defined as in Section 2 with
$\lambda = 1/4$, and let $X(\cdot )$ be defined by (\ref{4.5e}). Let
$u_a$, $A_a$, $R_1$ and $R_2$ be sufficiently regular (for the
following estimates to make sense) and satisfy the estimates
\beq
\label{4.6e}
\parallel u_a(t); W_{\infty}^2 \parallel \ \vee \ \parallel \partial_t u_a(t) \parallel_{\infty} \ \leq c\ t^{-3/2} \ ,
\eeq
\beq
\label{4.7e}
\parallel \partial_t^j A_a(t) \parallel_{\infty} \ \leq a \ t^{-1} \qquad \hbox{\it for $j = 0, 1$}\ ,
\eeq  
\beq
\label{4.8e}
\parallel \partial_t^j R_1; L^1([t , \infty ), L^2)\parallel \ \leq r_1 \ h(t) \qquad \hbox{\it for $j = 0, 1$}\ ,
\eeq 
\beq
\label{4.9e}
\parallel R_1; L^{8/3}([t , \infty ), L^4)\parallel \ \leq r_1 \ t^{-\eta}\ h(t) \qquad \hbox{\it for some $\eta \geq 0$}\ ,
\eeq 
\beq
\label{4.10e}
\parallel \omega^{-1}R_2; L^1([t , \infty ), H^1)\parallel \ \leq r_2 \ h(t) \ ,
\eeq

\noi for some constants $c$, $a$, $r_1$ and $r_2$ and for all $t \geq
1$. Then there exists $T$, $1 \leq T < \infty$, and there exists a
unique solution $(v, B)$ of the system (\ref{4.2e}) in $X([T, \infty
))$. If in addition
\beq
\label{4.11e}
\parallel \omega^{-1}R_2; L^1([t , \infty ), L^2)\parallel \ \leq r_2 \ t^{-1/2}\ h(t) 
\eeq

\noi for all $t \geq T$, then $B$ satisfies the estimate
\beq
\label{4.12e}
\parallel B(t);H^1\parallel\ \vee \ \parallel \omega^{-1}\partial_t B(t); H^1)\parallel \ \leq C\left ( t^{-1/2} + t^{1/4}  \ h(t)\right ) \ h(t) 
\eeq

\noi for some constant $C$ and for all $t \geq T$.}\\

The proof follows again those of Propositions 2.1 and 2.4. One starts
from $(v, B) \in X([T, \infty))$ for some $T$, $1 \leq T < \infty$, so
that $(v, B)$ satisfies 
\beq
\label{4.13e}
\left \{ \begin{array}{l}
\parallel v(t) \parallel_2 \ \leq N_0\ h(t)\\
\\
\parallel v ; L^{8/3}(J, L^4)\parallel \ \leq \  N_1\ h(t)\\ 
\\
\parallel B(t);H^1 \parallel \ \vee \ \parallel \partial_t B(t)\parallel_2\ \leq \  N_2\ h(t) \\
\\
\parallel \partial_t v(t) \parallel_2 \ \leq \  N_3\ h(t)\\
\\
\parallel \partial_t v ; L^{8/3}(J, L^4)\parallel \ \leq \  N_4\ h(t)\\ 
\\
\parallel \Delta v(t) \parallel_2\ \leq N_5\ h(t) \\
\\
\parallel \Delta v; L^{8/3}(J, L^4) \parallel\ \leq N_6\ h(t) \\
\end{array}\right . \eeq

\noi for some constants $N_i$, $0 \leq i \leq 6$ and for all $t \geq
T$, with $J = [t, \infty )$. For each such $(v, B)$ one constructs a
solution $(v',B')$ of the system (\ref{4.4e}) in $X([T, \infty ))$ and
one shows that the map $\phi : (v, B) \to (v', B')$ thereby defined is
a contraction on the subset ${\cal R}$ of $X([T, \infty))$ defined by
(\ref{4.13e}) for suitably chosen $N_i$ and for sufficiently large $T$.
The fact that the (Z)$_3$ system is short range manifests itself at
this stage by the absence of smallness conditions on $u_a$.
Technically, this follows from the fact that the nontriangular linear
terms in the system of equations for the $N_i$ which ensures the
contraction contain decaying powers of $T$ and can be made small by
taking $T$ large.\par

We now turn to Step 2, namely to the construction of $(u_a, A_a)$
satisfying the assumptions of Proposition 4.1. The (Z)$_3$ system is
short range and one can take for $(u_a, A_a)$ a pair of solutions
$(u_0, A_0)$ of the underlying free linear system, namely
\beq
\label{4.14e}
u_0(t) = U(t) u_+
\eeq

\noi and $A_0$ defined by (\ref{2.24e}). The remainders become  
\beq
\label{4.15e}
\left \{ \begin{array}{l}
R_1 = - A_0\ u_0\\
\\
R_2 = - \Delta |u_0|^2
\end{array}\right . \eeq

\noi and the final result can be stated as follows.\\

\noi {\bf Proposition 4.2.} {\it Let $h(t) = t^{-1/2}$ and let $X(\cdot
)$ be defined by (\ref{4.5e}). Let $u_+ \in H^2 \cap W_1^2$, let $A_+,
\omega^{-1} \dot{A}_+ \in H^1$ and $\nabla^2 A_+, \nabla \dot{A}_+ \in
W_1^1$. Let $(u_0, A_0)$ be defined by (\ref{4.14e}) (\ref{2.24e}).
Then there exists $T$, $1 \leq T < \infty$, and there exists a unique
solution $(u,A)$ of the (Z)$_3$ system (\ref{4.1e}) such that $(v,B) \equiv
(u-u_0, A- A_0) \in X([T, \infty)$. If in addition $u_+ \in H^{0,2}$,
then $B$ satisfies the estimate
\beq
\label{4.16e}
\parallel B(t);H^1\parallel\ \vee \ \parallel \omega^{-1}\partial_t B(t); H^1)\parallel \ \leq C\ t^{-3/4}  
\eeq
 
 \noi for some constant $C$ and for all $t \geq T$.}\\
 
The result follows from Proposition 4.1, from the dispersive estimate
(\ref{2.37e}), from Lemma 2.3 and Sobolev inequalities. In particular $$\parallel R_1 \parallel_2 \ \leq \ \parallel A_0 \parallel_2 \ \parallel u_0\parallel_{\infty} \ \leq C\ t^{-3/2}$$ 

\noi and similar estimates lead to the decay $h(t) = t^{-1/2}$, while 
$$\parallel R_1 \parallel_4 \ \leq \ \parallel A_0 \parallel_4 \ \parallel u_0\parallel_{\infty} \ \leq C\ t^{-2}$$ 

\noi ensures the condition (\ref{4.9e}) with $\eta = 9/8$. The last
statement of Proposition 4.2 requires in addition the following lemma,
which we state in dimensions $n = 2, 3$, and which follows immediately
from the factorisation (\ref{2.26e}) of $U(t)$.\\

\noi {\bf Lemma 4.1.} {\it  Let $n = 2$ or $3$. Let $u_+ \in
H^{0,2}(\subset L^1)$ and let $u_0 = U(t) u_+$. Then the following estimates
hold~:}
\beq
\label{4.17e}
\parallel \nabla |u_0|^2\parallel_2 \ \leq 2 (2 \pi |t|)^{-n/2} \ t^{-1} \parallel u_+\parallel_1 \ \parallel xu_+\parallel_2 \ ,
\eeq
\beq
\label{4.18e}
\parallel \Delta |u_0|^2\parallel_2 \ \leq 4 (2 \pi t)^{-n/2} \ t^{-2} \parallel u_+\parallel_1 \ \parallel x^2u_+\parallel_2 \ .
\eeq
\vskip 5 truemm

In the same way as for the (WS)$_3$ system, by using a more accurate
asymptotic form for $u_a$, one can obtain a stronger asymptotic
convergence in time of the solution on a smaller subspace of asymptotic
states \cite{34r}. Thus we choose $(u_a, A_a)=$\break\noindent $((1 + f)u_0, A_0)$ with
$(u_0, A_0)$ defined by (\ref{4.14e}) (\ref{2.24e}) and $f$ by
(\ref{2.49e}), namely $f = 2 \Delta^{-1}A_0$. Using the operators $J$
and $P$ defined by (\ref{2.50e}), we rewrite $R_1$ as
\bea
\label{4.19e}
R_1 &=& \left ( i \partial_t + (1/2) \Delta - A_0\right ) (1 + f) u_0\nn \\
&=& - f\ A_0\ u_0 - it^{-1} (\nabla f)\cdot Ju_0 + i t^{-1} (Pf) u_0
\eea

\noi while $R_2$ now becomes
\beq
\label{4.20e}
R_2 = - \Delta (1 + f)^2|u_0|^2\ .
\eeq

\noi The new remainders are easily estimated through the following lemmas.\\

\noi {\bf Lemma 4.2.} {\it Let $u_+ \in W_1^2$, $xu_+ \in W_1^2$, and let $A_0$ and $f$ satisfy 
\beq
\label{4.21e}
\parallel \partial_t^j \nabla^k A_0\parallel_{\infty} \ \leq a\ t^{-1} 
\eeq
\beq
\label{4.22e}
\parallel \partial_t^j \nabla^k f;H^1\parallel\ \vee \  \parallel \partial_t^j \nabla^k Pf\parallel_2\ \leq C 
\eeq

\noi for $0 \leq j + k \leq 1$ and for all $t \geq 1$. Then the following estimates hold~:
\beq
\label{4.23e}
\parallel \partial_t^j \nabla^k R_1\parallel_2 \ \leq C\ t^{-5/2} 
\eeq

\noi for some constant $C$, for $0 \leq j + k \leq 1$ and for all $t \geq 1$.}\\

\noi {\bf Lemma 4.3.} {\it Let $u_+ \in W_1^2 \cap H^{0,2}$ and let $f$ satisfy 
\beq
\label{4.24e}
\parallel \nabla f(t)\parallel_2\ \vee \  \parallel \Delta f(t)\parallel_2\ \vee \ \parallel f(t) \parallel_{\infty}\ \leq C 
\eeq

\noi for all $t \geq 1$. Then the following estimates hold~: 
\beq
\label{4.25e}
\parallel \omega^{-1}R_2 \parallel_2 \ \leq C\ t^{-5/2} \ ,
\eeq
\beq
\label{4.26e}
\parallel R_2 \parallel_2 \ \leq C\ t^{-3} 
\eeq

\noi for some constant $C$ and for all $t \geq 1$.}\\

\noi Lemma 4.3 follows readily from Lemma 4.1. In practice, the bound on
$\parallel f\parallel_{\infty}$ in (\ref{4.24e}) will follow from the
Sobolev inequality
$$\parallel f \parallel_{\infty} \ \leq \ C \left ( \parallel \nabla f \parallel_2 \ \parallel \Delta f \parallel_2 \right )^{1/2}$$

\noi for $f$ tending to zero at infinity in some weak sense.\par

We can now state the final result with improved asymptotics.\\

\noi {\bf Proposition 4.3.} {\it Let $h(t) = t^{-3/2}$ and let $X(\cdot
)$ be defined by (\ref{4.5e}). Let $u_+ \in H^2 \cap H^{0,2} \cap
W_1^2$ with $xu_+ \in W_1^2$. Let $(A_+,\dot{A}_+)$ satisfy
\bea
\label{4.27e}
&&A_+, \omega^{-1} \dot{A}_+ \in \dot{H}^{-2} \cap H^1\quad , \quad \nabla^2A_+, \nabla \dot{A}_+ \in W_1^1 \ ,\nn\\
&&x\cdot \nabla A_+, \omega^{-1} x \cdot \nabla \dot{A}_+ \in \dot{H}^{-2} \cap \dot{H}^{-1}\ .
\eea

\noi Let $(u_0, A_0)$ be defined by (\ref{4.14e}) (\ref{2.24e}) and let
$u_a = (1 + f)u_0$ with $f = 2 \Delta^{-1}A_0$. Then :\par

(1) There exists $T$,
$1 \leq T < \infty$ and there exists a unique solution $(u, A)$ of the
(Z)$_3$ system (\ref{4.1e}) such that $(v, B) \equiv (u - u_a, A - A_0) \in
X([T, \infty ))$.\par

(2) Assume in addition that $\nabla^2 \omega^{-2} A_+$, $\nabla \omega^{-2} \dot{A}_+ \in W_1^1$. Then there exists $T$,
$1 \leq T < \infty$ and there exists a unique solution $(u, A)$ of the
(Z)$_3$ system (\ref{4.1e}) such that $(u - u_0, A - A_0) \in
X([T, \infty ))$. One can take the same $T$ and the solution $(u, A)$ is the same as in Part (1).}\\

As in Proposition 2.3, Part (2) follows immediately from Part (1) and
from the fact that $(fu_0, 0) \in X([T, \infty))$ for any $T \geq 1$.
The additional assumption on $(A_+, \dot{A}_+)$ has been made to ensure
that property, which does not seem to follow immediately from the
previous assumptions. It is stronger than needed. Together with the
previous assumptions and through Lemma 2.3, it implies that $f$,
$\nabla f$ and $\partial_tf$ also satisfy the decay estimates stated
for $A_0$ in (\ref{2.41e}) with $\ell = 0,1$. Only the special case $r
= 4$ is used for the present purpose. More economical assumptions could
be made instead by using suitable Besov spaces.\par

Note also that the assumption (\ref{4.27e}) on
$(A_+,\dot{A}_+)$ is the same as the assumption (\ref{2.69e}) of
Proposition 2.6 relative to the (WS)$_3$ system, a reflection of the
fact that we are treating $A_0$ in exactly the same way in both
cases.\par

We now turn to the case of space dimension $n = 2$, where the situation
is much less satisfactory. The free part of the asymptotic field is
estimated at best as
\beq
\label{4.28e}
\parallel A_0(t) \parallel_{\infty} \ \leq C\ t^{-1/2} 
\eeq

\noi and we are unable to handle such a slow decay in Step 1, so that
the final result will eventually be restricted to the special case of
zero asymptotic state $(A_+,\dot{A}_+)$ for $A$. On the other hand, in
a suitable limit, the Zakharov system formally yields the cubic NLS
equation, which is short range for $n=2$, and one might naively expect
a similar situation for the (Z)$_2$ system, allowing for a treatment of
that system without a smallness condition on $u$. This turns out not to
be the case, and the (Z)$_2$ system does actually require such a
smallness condition at the level of Step 1.  The treatment of that step
is very similar to the case of (Z)$_3$. The relevant space $X(\cdot )$
is essentially the same up to obvious changes, namely

\bea
\label{4.29e}
&&X(I) = \Big \{ (v, B):(v,B)\in {\cal C}(I, H^2\oplus H^1) \cap {\cal C}^1(I, L^2 \oplus L^2), \nn \\
&&\parallel (v, B);X(I)\parallel \ \equiv \ \mathrel{\mathop {\rm Sup}_{t \in I }}\ h(t)^{-1} 
\Big ( \parallel v(t);H^2\parallel \ + \ \parallel \partial_t v(t)\parallel_2\nn \\
&&+\ \parallel v; L^{4}(J,W_4^2)\parallel \ + \ \parallel \partial_t v; L^{4}(J, L^4)\parallel \nn\\
&& + \ \parallel B(t); H^1\parallel\ + \ \parallel \partial_t B(t) \parallel_2 \Big ) < \infty \Big \} 
\eea

\noi where $J = [t, \infty ) \cap I$, and the result can be stated as follows.\\

\noi {\bf Proposition 4.4} {\it Let $h$ be defined as in Section 2 with
$\lambda = 1/2$ and let $X(\cdot )$ be defined by (\ref{4.29e}). Let
$u_a$, $A_a$, $R_1$ and $R_2$ be sufficiently regular and satisfy the estimates
\beq
\label{4.30e}
\parallel u_a(t); W_{\infty}^2 \parallel \ \vee \ \parallel \partial_t u_a(t) \parallel_{\infty} \ \leq c\ t^{-1} \ ,
\eeq
\beq
\label{4.31e}
\parallel \partial_t^j A_a(t) \parallel_{\infty} \ \leq a \ t^{-1-j\theta} \qquad \hbox{\it for some $\theta >0$ and for $j = 0, 1$}\ ,
\eeq  
\beq
\label{4.32e}
\parallel \partial_t^j R_1; L^1([t , \infty ), L^2)\parallel \ \leq r_1 \ h(t) \qquad \hbox{\it for $j = 0, 1$}\ ,
\eeq 
\beq
\label{4.33e}
\parallel R_1; L^{4}([t , \infty ), L^4)\parallel \ \leq r_1 \ t^{-\eta}\ h(t) \qquad \hbox{\it for some $\eta \geq 0$}\ ,
\eeq 
\beq
\label{4.34e}
\parallel \omega^{-1}R_2; L^1([t , \infty ), H^1)\parallel \ \leq r_2 \ h(t) 
\eeq

\noi for some constants $c$, $a$, $r_1$ and $r_2$ with $c$ sufficiently small and for all $t \geq
1$. Then there exists $T$, $1 \leq T < \infty$, and there exists a
unique solution $(v, B)$ of the system (\ref{4.2e}) in $X([T, \infty
))$.}\\

The proof is a minor variation of that of Proposition 4.1. The
assumption (\ref{4.31e}) is rather arbitrary. It is too strong to
accomodate a non zero $A_0$ (see (\ref{4.28e})) but weaker by one power
of $t$ than the condition that would be satisfied by an $A_1$ of the
type (\ref{2.25e}) devised to ensure $R_2 = 0$. It has been chosen so
as to ensure that the proof of the proposition proceeds smoothly. \par

The final result can be stated as follows.\\
 
\noi {\bf Proposition 4.5.} {\it Let $h(t) = t^{-1}$ and let $X(\cdot
)$ be defined by (\ref{4.29e}). Let $u_+ \in H^2 \cap H^{0,2} \cap
W_1^2$ with $\parallel u_+ ; W_1^2\parallel$ sufficiently small and let $u_0(t) = U(t) u_+$.
 Then there exists $T$,
$1 \leq T < \infty$, and there exists a unique solution $(u, A)$ of the
(Z)$_2$ system (\ref{4.1e}) such that $(u - u_0, A ) \in
X([T, \infty ))$.}\\

The result follows readily from Proposition 4.4 with $A_a = 0$, from
the dispersive estimate (\ref{2.37e}) and from Lemma 4.1.

\mysection{The Klein-Gordon-Schr\"odinger system (KGS)$_{2}$}
\hspace*{\parindent} In this section we review the main results available on the local
Cauchy problem at infinity for the (KGS)$_2$ system
\beq
\label{5.1e}
\left \{ \begin{array}{l} i \partial_t u = - (1/2) \Delta u + A u\\ \\(\sq + 1) A = - |u|^2\ . \end{array} \right .
\eeq

That system is known to be globally well posed in the energy space
\cite{1r} \cite{4r}. We shall first present the results concerning Step
1, which are easily obtained by minor variations from the corresponding
results for the (WS)$_3$ and (Z)$_2$ systems, especially from
Propositions 2.1, 2.4 and 4.4. We shall then outline the construction
of the asymptotic form $(u_a, A_a)$ given in \cite{29r}. That
construction is more delicate than in the (WS)$_3$ case. We  shall then
state the final results in a partly qualitative way, but we shall refrain from a completely formal
statement. The reason is that Step 1 is treated in \cite{29r} by a
slightly different method and in spaces smaller (more regular) than
here, so that the assumptions made in \cite{29r} on the asymptotic
state are stronger than needed for a combination with the treatment of
Step 1 given here. We leave it as an open question to perform that
combination and in particular to determine the most economical and
natural assumptions on the asymptotic state needed for that
purpose.\par

We first consider Step 1. We follow again the sketch given in the
introduction. For a given asymptotic $(u_a, A_a)$, we look for $(u, A)$
in the form $(u, A) = (u_a + v, A_a + B)$. The auxiliary system
satisfied by $(v, B)$ is now
\beq
\label{5.2e}
\left \{ \begin{array}{l} i \partial_t v = - (1/2) \Delta v + A v + Bu_a - R_1 \\ \\(\sq + 1) B = - \left ( |v|^2 + 2\ {\rm Re}\ \overline{u}_a v\right ) - R_2 \end{array} \right .
\eeq

\noi where the remainders are defined by
\beq
\label{5.3e}
\left \{ \begin{array}{l} R_1 = i \partial_t u_a + (1/2) \Delta u_a - A_a u_a \\ \\R_2 = (\sq + 1) A_a +|u_a|^2 \ . \end{array} \right .
\eeq

\noi We consider also the partly linearized system for functions $(v', B')$  
\beq
\label{5.4e}
\left \{ \begin{array}{l} i \partial_t v' = - (1/2) \Delta v' + A v' + Bu_a - R_1 \\ \\(\sq + 1) B' = - \left ( |v|^2 + 2\ {\rm Re}\ \overline{u}_a v\right ) - R_2 \ .\end{array} \right .
\eeq

\noi Step 1 of the method consists again in solving the system
(\ref{5.2e}) for $(v, B)$ with $(v, B)$ tending to zero at infinity
under assumptions on $(u_a, A_a)$ of a general nature. As in the case of
the (Z)$_n$ system, but now for a different reason, namely because the
KG energy controls the $L^2$ norm, it is easy to estimate the $L^2$
norm of $B$ (through the energy), so that the Strichartz inequalities for
the KG equation do not seem to be useful here. On the other hand, as
in the case of the (WS)$_3$ system and in contrast with the case of
(Z)$_n$, the low level of local singularity of (KGS)$_2$ allows for a
treatment of Step 1 in low level spaces, and in particular with $k = 0$. The
largest suitable space $X(\cdot )$ can be defined as follows 
\bea
\label{5.5e}
&&X(I) = \Big \{ (v, B):(v,B)\in {\cal C}(I, L^2\oplus H^1), B \in  {\cal C}^1(I, L^2), \nn \\
&&\parallel (v, B);X(I)\parallel \ \equiv \ \mathrel{\mathop {\rm Sup}_{t \in I }}\ h(t)^{-1} 
\Big ( \parallel v(t)\parallel_2 \ + \ \parallel v; L^{4}(J,L^4)\parallel \nn \\
&& + \ \parallel B(t); H^1\parallel\ + \ \parallel \partial_t B(t) \parallel_2 \Big ) < \infty \Big \} 
\eea

\noi where $J = [t, \infty ) \cap I$ and $h$ is defined as in Section
2. That space has $(k, \ell ) = (0,1)$. The definition includes the
$L^2$ norm and a Strichartz norm of level $k = 0$ for $v$, and the
energy norm for $B$, which has level $\ell = 1$. \par

We now state the result concerning Step 1, which is a minor variation of Proposition 2.1\\

\noi {\bf Proposition 5.1} {\it Let $h$ be defined as in Section 2 with
$\lambda = 1/2$ and let $X(\cdot )$ be defined by (\ref{5.5e}). Let
$u_a$, $A_a$, $R_1$ and $R_2$ be sufficiently regular (for the following estimates to
make sense) and satisfy the estimates
\beq
\label{5.6e}
\parallel u_a(t) \parallel_{\infty}\ \leq c_0 \ t^{-1}\ ,
\eeq
\beq
\label{5.7e}
\parallel A_a(t) \parallel_{\infty} \ \leq a\ t^{-1}\ ,
\eeq
\beq
\label{5.8e}
\parallel R_1;L^1([t, \infty ), L^2) \parallel\ \leq r_1\ h(t)\ ,
\eeq
\beq
\label{5.9e}
\parallel R_2;L^{1}([t, \infty ), L^{2}) \parallel\ \leq r_2\ h(t)
\eeq

\noi for some constants $c_0$, $a$, $r_1$ and $r_2$ with $c_0$
sufficiently small and for all $t \geq 1$. Then there exists $T$, $1
\leq T < \infty$ and there exists a unique solution $(v, B)$ of the
system (\ref{5.2e}) in $X([T, \infty ))$.} \\

The proof follows closely that of Proposition 2.1. One starts from $(v, B) \in X([T, \infty))$ for some $T$, $1 \leq T < \infty$, so that $(v,B)$ satisfies 
\beq
\label{5.10e}
\left \{ \begin{array}{l} \parallel v(t) \parallel_2 \ \leq N_0\ h(t)\\
\\
\parallel v ; L^{4}([t, \infty), L^4) \parallel \ \leq N_1\ h(t)\\
\\
\parallel B(t) ; H^{1}\parallel \ \vee \ \parallel \partial_t B(t)\parallel_2\ \leq N_2\ h(t)
\end{array}\right .\eeq
 
\noi for some constants $N_i$, $0 \leq i \leq 2$ and for all $t \geq
T$. For each such $(v, B)$ one constructs a solution $(v', B')$ of the
system (\ref{5.4e}) in $X([T, \infty ))$ and one shows that the map
$\phi : (v, B) \to (v', B')$ thereby defined is a contraction on the
subset ${\cal R}$ of $X([T, \infty ))$ defined by (\ref{5.10e}) for
suitably chosen $N_i$ and for sufficiently large $T$. The seminorms
$N'_i$ of $(v', B')$ corresponding to (\ref{5.10e}) are estimated by
the use of (\ref{2.3e}) (\ref{2.7e}) (\ref{2.9e}) as 
\beq
\label{5.11e}
\left \{ \begin{array}{l} N'_0 \leq 2c_0 \ N_2 + r_1  \\
\\
N'_1 \leq C_1 \left ( c_0 \ N_2 + r_1\right ) ( 1 + a ) \\
\\
N'_2 \leq 4 c_0 \ N_0 + r_2 + 4N_1^2\ \overline{h}(T) 
\end{array} \right .
\eeq

\noi for $T$ sufficiently large to ensure that $C_1 N_2 \overline{h}(T)
\leq 1$, for some absolute constant $C_1$. Sufficient conditions to
ensure the stability of ${\cal R}$ under $\phi$ are obtained in analogy
with (\ref{2.21e}) by imposing
\beq
\label{5.12e}
\left \{ \begin{array}{l} N_0 = 2c_0 \ N_2 + r_1  \\
\\
N_1 = C_1 \left ( c_0 \ N_2 + r_1\right ) ( 1 + a ) \\
\\
N_2 = 4 c_0 \ N_0 + r_2 + 1 
\end{array} \right .
\eeq

\noi which is possible under the smallness condition $8c_0^2 < 1$ and
by taking $T$ sufficiently large so that in addition
$4N_1^2\overline{h}(T)  \leq 1$. The rest of the proof is a minor
variation of that of Proposition 2.1.\\

As in the case of the (WS)$_3$ system, there is no difficulty to
implement Step 1 at higher levels of regularity. As in the case of
(WS)$_3$, the time decay remains the same, namely $\lambda = 1/2$, and
no additional smallness condition appears beyond the previous one of
$c_0$. Again two theories are of special interest. \par

(1) The theory at the level $(k , \ell ) = (1, 1)$ of the energy. The
appropriate function space is now 
\bea
\label{5.13e}
&&X_1(I) = \Big \{ (v, B):(v,B)\in {\cal C}(I, H^1\oplus H^1), B \in  {\cal C}^1(I, L^2), \nn \\
&&\parallel (v, B);X_1(I)\parallel \ \equiv \ \mathrel{\mathop {\rm Sup}_{t \in I }}\ h(t)^{-1} 
\Big ( \parallel v(t);H^1\parallel \ + \ \parallel v; L^{4}(J,W^1_4)\parallel \nn \\
&& + \ \parallel B(t); H^1\parallel\ + \ \parallel \partial_t B(t) \parallel_2 \Big ) < \infty \Big \} 
\eea

\noi where $J = [t, \infty) \cap I$. It differs from the previous $X(I)$ by the inclusion of the $L^2$
norm and of a Strichartz norm of $v$ at the level $k=1$. The additional
estimates needed in that theory are obtained from (\ref{2.5e}) with
$\partial = \nabla$ and from (\ref{2.9e}) applied to $\nabla v$. \par

(2) The theory at the level $(k, \ell ) = (2, 1)$ because that is the
lowest convenient level appropriate for the (Z)$_2$ system. Actually
one can make the same choice of $X(I)$ as in the latter case, namely
\bea
\label{5.14e}
&&X_2(I) = \Big \{ (v, B):(v,B)\in {\cal C}(I, H^2\oplus H^1) \cap {\cal C}^1(I, L^2 \oplus L^2), \nn \\
&&\parallel (v, B);X_2(I)\parallel \ \equiv \ \mathrel{\mathop {\rm Sup}_{t \in I }}\ h(t)^{-1} 
\Big ( \parallel v(t);H^2\parallel \ + \ \parallel \partial_t v(t)\parallel_2\nn \\
&&+\ \parallel v; L^{4}(J,W_4^2)\parallel \ + \ \parallel \partial_t v; L^{4}(J, L^4)\parallel \nn\\
&& + \ \parallel B(t); H^1\parallel\ + \ \parallel \partial_t B(t) \parallel_2 \Big ) < \infty \Big \} 
\eea

\noi which coincides with (\ref{4.29e}). We state the result concerning
that theory in order to allow for comparison both with the case of the
(Z)$_2$ system treated in Section 4 and with the treatment of Step 1 in
\cite{29r} which uses a space of level $(k, \ell ) = (2,2)$.\\

\noi {\bf Proposition 5.2.} {\it Let $h$ be defined as in Section 2 with
$\lambda = 1/2$ and let $X_2(\cdot )$ be defined by (\ref{5.14e}). Let
$u_a$, $A_a$, $R_1$ and $R_2$ be sufficiently regular (for the following estimates to
make sense) and satisfy the estimates
\beq
\label{5.15e}
\parallel \partial_t^j u_a(t) \parallel_{\infty}\ \leq c \ t^{-1}\qquad \hbox{\it for $j = 0, 1$}
\eeq

\noi and in particular (\ref{5.6e}),
\beq
\label{5.16e}
\parallel \partial_t^j A_a(t) \parallel_{\infty} \ \leq a\ t^{-1}\qquad \hbox{\it for $j = 0, 1$}\ ,
\eeq
\beq
\label{5.17e}
\parallel \partial_t^j R_1;L^1([t, \infty ), L^2) \parallel\ \leq r_1\ h(t)\qquad \hbox{\it for $j = 0, 1$}\ ,
\eeq
\beq
\label{5.18e}
\parallel R_1;L^{4}([t, \infty ), L^{4}) \parallel\ \leq r_1\ t^{-\eta }\ h(t) \qquad \hbox{\it for some $\eta \geq 0$}\ ,
\eeq
$$\parallel R_2;L^{1}([t, \infty ), L^{2}) \parallel\ \leq r_2\ h(t) \eqno(5.9)\equiv(5.19)$$

\noi for some constants $c$, $c_0$, $a$, $r_1$ and $r_2$ with $c_0$ 
sufficiently small and for all $t \geq 1$. Then there exists $T$, $1
\leq T < \infty$ and there exists a unique solution $(v, B)$ of the
system (\ref{5.2e}) in $X_2([T, \infty ))$.} \\  

The proof is a combination of those of Propositions 5.1 and 4.4. One
starts from $(v, B) \in X_2 ([T, \infty ))$ for some $T$, $1 \leq T <
\infty$, so that $(v, B)$ satisfies (\ref{5.10e}) and in addition
$$\left \{ \begin{array}{l} \parallel \partial_t v(t) \parallel_2 \ \leq N_3\ h(t)\\
\\
\parallel \partial_t v ; L^{4}([t, \infty ), L^4) \parallel \ \leq N_4\ h(t)\\
\\
\parallel \Delta v(t) \parallel_2 \ \leq N_5\ h(t)\\
\\
\parallel \Delta v ; L^{4}([t, \infty ), L^4)\parallel \ \leq N_6\ h(t)
\end{array}\right . \eqno(5.20)$$

\noi for some constants $N_i$, $0 \leq i \leq 6$ and for all $t \geq
T$. For each such $(v, B)$ one constructs a solution $(v', B')$ of the
system (\ref{5.4e}) in $X([T, \infty ))$. Now the additional estimates
of the norms of $v'$ corresponding to (5.20) are identical with
those that occur in the proof of Proposition 4.4 with $\theta = 0$,
since one is estimating the same norms of $v'$ from the same
Schr\"odinger equation and with the same information on $A_a$, $B$ and
$R_1$. The rest of the proof proceeds as before.\\

We now turn to Step 2, namely the construction of $(u_a, A_a)$
satisfying the assumptions needed for Step 1, and we describe the
choice made in \cite{29r}. As in the (WS)$_3$ case, one sees
immediately that taking for $(u_a, A_a)$ a pair of solutions of the
free equations is inadequate. Furthermore introducing a phase in $u_a$
does not seem to be useful at this stage. One takes 
$$u_a = u_1 + u_2 \equiv (1 + f) u_1 \ , \eqno(5.21)$$

\noi where $f$ is a complex valued function to be chosen later and
$$u_1 = MD\widehat{u}_+ \ . \eqno(5.22)$$

\noi The difference between $u_1$ and $u_0 = U(t)u_+$ can be written as
$$u_0 - u_1 = U(t) \left ( 1 - \overline{M}(t)\right ) u_+ \eqno(5.23)$$

\noi and is easily controlled under suitable assumptions on $u_+$. One then takes
$$A_a = A_0 + A_1 \eqno(5.24)$$
$$A_0 = \cos \omega_1 t \ A_+ + \omega_1^{-1} \sin \omega_1 t\ \dot{A}_+\eqno(5.25)$$

\noi with $\omega_1 = (1 - \Delta )^{1/2}$, so that $(\sq + 1)A_0 = 0$, and one takes for $A_1$ the solution of the equation 
$$(\sq + 1) A_1 = - |u_1]^2 \eqno(5.26)$$

\noi which vanishes at infinity. Using an integration by parts in time, one can write $A_1$ as 
$$A_1(t)  = \int_t^{\infty} dt ' \left ( 1 - \cos \left ( \omega_1 (t'-t)\right ) \right ) \omega_1^{-2} \partial_t |u_1|^2(t') \ .\eqno(5.27)$$

\noi This formula is justified and a good control of $A_1$ is provided
by the following lemma, which is a variant of Lemma 3.2 in \cite{29r}
or Lemma 2.4 in \cite{28r}.\\

\noi {\bf Lemma 5.1}. {\it Let $g \in {\cal C}^1([1, \infty ), L^2)$
with $\partial_t g \in L^1([1, \infty ), L^2)$ so that $g(t)$ has an
$L^2$ limit $g(\infty )$ as $t \to \infty$ and let $g( \infty ) = 0$.
Then there exists a unique solution $A$ of the equation
$$(\sq + 1 )A = g$$

\noi such that $(A, \partial_t A) \in {\cal C}([1, \infty ), H^1 \oplus
L^2)$ and that $\parallel \omega_1 A \parallel_2  \vee \parallel
\partial_t A \parallel_2 \ \to 0$ as $t \to \infty$. Furthermore $(A,
\partial_t A) \in {\cal C} ([1, \infty ), H^2 \oplus H^1)$ and $A$
satisfies the estimates

$$\left \{ \begin{array}{l} \parallel \omega_1^2 A(t) \parallel_2 \ \leq \ 2\parallel \partial_t g ; L^{1}([t, \infty ), L^2) \parallel \\
\\
\parallel \omega_1 \partial_t A(t) \parallel_2 \ \leq \ \parallel \partial_t g ; L^{1}([t, \infty ), L^2) \parallel 
\end{array}\right . \eqno(5.28)$$

\noi for all $t \geq 1$.}\\

The proof uses an integration by parts in time and a limiting procedure. It yields in particular the representations
$$\left \{ \begin{array}{l} \omega_1^2 A(t) = - \displaystyle{\int_t^{\infty}} dt' \left ( 1 - \cos \left ( \omega_1(t'-t)\right ) \right ) \partial_t g (t')\\
\\
\omega_1 \partial_t A(t) =  \displaystyle{\int_t^{\infty}} dt' \sin \left ( \omega_1(t'-t)\right ) \partial_t g (t')
\end{array}\right . \eqno(5.29)$$

\noi of which (5.27) is a special case and from which the estimates (5.28) follow immediately. In the case of (5.26), one has
$$g(t) = - |u_1(t)|^2 = - t^{-2} D_0 (t) |\widehat{u}_+|^2 \eqno(5.30)$$

\noi so that by the commutation rule
$$\partial_t D_0(t) = t^{-1} D_0(t) \left ( t \partial_t - x \cdot \nabla \right ) \eqno(5.31)$$

\noi one obtains
$$\partial_t g(t) = - \partial_t |u_1 (t)|^2 = t^{-3} D_0(t) \left ( 2|\widehat{u}_+|^2 +  x \cdot \nabla |\widehat{u}_+|^2 \right ) \eqno(5.32)$$

\noi and therefore by (5.28)
$$\parallel \omega_1^2 A_1(t) \parallel_2 \ + \ \parallel \omega_1 \partial_t A_1(t) \parallel_2\ \leq 3t^{-1}\left ( 2\parallel \widehat{u}_+\parallel_4^2 \ + \ \parallel x \cdot \nabla |\widehat{u}_+|^2 \parallel_2 \right )\ . \eqno(5.33)$$

\noi In particular $A_1$ exhibits a $t^{-1}$ decay in norms for which $A_0$ has no decay. \par

We now turn to the remainders. Substituting (5.21) (5.24) into their definition (\ref{5.3e})
and using again $J$ and $P$ defined by (\ref{2.50e}), we obtain
(compare with (\ref{2.48e}) (\ref{2.51e}))
$$R_1 = (1 + f) \left ( i \partial_t + (1/2) \Delta -A_1\right ) u_1 - f A_0 u_1 - i t^{-1} \nabla f \cdot Ju_1$$
$$+ \left ( (1/2) \Delta f + i t^{-1} P f - A_0 \right ) u_1 \ , \eqno(5.34)$$
$$R_2 = \left ( |f|^2 + 2\ {\rm Re}\ f\right ) |u_1|^2 \ . \eqno(5.35)$$

\noi The first three terms in the RHS of (5.34) and the RHS of (5.35) are expected and will turn out to be
$O(t^{-2})$ in the relevant norms, thereby fulfilling the required
assumptions of Propositions 5.1 and 5.2 with $h(t) = t^{-1}$. On the other hand the last term in (5.34) has only $t^{-1}$ decay
for general $f$ and one has to choose $f$ in order to improve that
decay to $t^{-2}$. However this is more delicate than in the (WS)$_3$
case, because if $f$ is a solution of the free KG equation, in general
$Pf$ is not, so that $Pf$ has to be taken into account in order to
perform the cancellation. This is done by using the explicit asymptotic
form of solutions of the free KG equation \cite{20newref}. It is
convenient to decompose $A_0$ into positive and negative frequency
parts, namely
$$A_0 = A_{0+} + A_{0-} = 2\ {\rm Re}\ A_{0+} \eqno(5.36)$$

\noi where
$$ A_{0\pm} = (1/2) \left ( A_0 \mp i \omega_1^{-1} \partial_t A_0\right ) \eqno(5.37)$$

\noi and $A_{0-} = \overline{A}_{0+}$ since $A_0$ is real. The equation satisfied by $A_{0\pm}$ is 
$$\left ( i \partial_t \pm \omega_1\right ) A_{0\pm} = 0 \eqno(5.38)$$

\noi with the initial condition
$$A_{0\pm} (0) = (1/2) \left ( A_+ \mp i \omega^{-1} \dot{A}_+\right ) \equiv A_{+\pm}\ . \eqno(5.39)$$

\noi The asymptotic form of $A_{0\pm}$ is then (\cite{20newref}), Theorem 7.2.5)
$$A_{0\pm} \sim \widetilde{A}_{0\pm} = \pm i (2 \pi t)^{-1} (t^2/\rho^2) \widehat{A}_{+\pm} (\mp x/\rho) \exp (\pm i \rho ) \chi (|x| < t)\eqno(5.40)$$

\noi where $\rho = (t^2 -
x^2)^{1/2}$ and $\chi (|x|< t)$ denotes the characteristic function
of the set $\{(x, t): |x| < t \}$. Corresponding to the
decomposition (5.36) of $A_0$, we look for $f$ in the form 
$$\left \{ \begin{array}{l} f = f_+ + f_-  \ , \\ \\
f_{\pm} = t^{-1} \ F_{\pm} (x/t) \exp (\pm i \rho ) \ .\end{array}\right . \eqno(5.41)$$

\noi We compute
$$(1/2) \Delta f_{\pm} + i t^{-1} Pf_{\pm} = t^{-1} \Big \{ - i t^{-1} F_{\pm}(x/t) + (2t^2)^{-1} (\Delta F_{\pm})(x/t)$$
$$\pm \ i t^{-1} (\nabla F_{\pm})(x/t) \cdot \nabla \rho - \left ( \pm \ t^{-1} P\rho + (1/2) |\nabla \rho|^2 \mp (i/2) \Delta \rho \right ) F_{\pm} (x/t)\Big \} \exp (\pm \ i \rho )\eqno(5.42)$$

\noi where we have used the fact that $P(F(x/t)) = 0$, 
$$\cdots = t^{-1} \Big \{ - i \left ( t^{-1} \pm \rho^{-1} \pm x^2(2\rho^3)^{-1} \right )F_{\pm} (x/t) + (2t^2)^{-1} (\Delta F_{\pm}) (x/t) $$
$$\mp \ i \rho^{-1} (x \cdot \nabla F_{\pm}) (x/t) - \left ( \pm \ \rho t^{-1}  + x^2 (2 \rho^2)^{-1} \right )F_{\pm}(x/t)\Big \} \exp (\pm \ i\rho ) \eqno(5.43)$$

\noi where we have used the fact that $P\rho = \rho$ and $\Delta
\rho = - (2/\rho + x^2/ \rho^3)$. Now all the terms in the bracket in
(5.43) are $O(t^{-1})$ in the sense that they become $Ct^{-1}$ for $x =
ct$, $|c| < 1$, except for the last term which is $O(1)$ in the same
sense. We now choose $F_{\pm}$ in such a way that this term cancels the
contribution of $\widetilde{A}_0$ in the last term of (5.34), namely 
$$\left ( \pm \ \rho  t^{-1} + x^2(2\rho^2)^{-1}\right ) F_{\pm}(x/t) = \pm \ (2 \pi i)^{-1} (t^2/\rho^2) \widehat{A}_{+\pm} (\mp \ x/\rho ) \chi (|x| < t)$$

\noi or equivalently
$$\left ( \pm \ (1 - x^2)^{3/2} + x^2/2 \right ) F_{\pm}(x) = \pm \ (2 \pi i)^{-1} \widehat{A}_{+\pm} \left ( \mp \ x(1- x^2)^{-1/2}\right )  \chi (|x| < 1)\ . \eqno(5.44)$$

\noi The functions $F_{\pm}$ defined by (5.44) can exhibit two kinds
of singularities. First, because of the occurrence of $(1 -
x^2)^{1/2}$ and of $\chi (|x| < 1)$, their derivatives could be
singular when $|x| \to 1$. However in that limit the argument of
$\widehat{A}_{+\pm}$ tends to infinity, and those singularities can be
controlled by assuming sufficient decay of $\widehat{A}_{+\pm}$ at
infinity. Second, $F_-$ contains the factor
$$g(x) = \left ( (1-x^2)^{3/2} - x^2/2\right )^{-1}$$

\noi which is singular on the cercle $C_r = \{x:|x|= r\}$ for a
suitable $r$, $0 < r < 1$. This singularity always appears in the
remainders in the form of combinations of the type $\partial^{\alpha}f
\partial^{\beta}u_1$ which generate combinations of the type
$\partial^{\alpha}F_- \partial^{\beta}\widehat{u}_+$ and it can
therefore be controlled by assuming that $\widehat{u}_+$ vanishes of
sufficient order on $C_r$. We refer to \cite{29r} for the estimates of
$u_a$ and of the remainders with the previous choice of $(u_a, A_a)$.
\par

Combining that choice with the treatment of
Step 1 provided by Propositions 5.1 and 5.2, one should obtain the
following final result, which we state in a semi qualitative way.\\

\noi {\bf Pseudoproposition 5.3.} {\it Let $h(t) = t^{-1}$ and let
$X(\cdot )$ be defined by (\ref{5.5e}). Let $(u_a, A_a)$ be defined by
(5.21) (5.22) (5.24) (5.25) (5.27) (5.41) (5.44) (5.39). Let $(u_+,
A_+, \dot{A}_+)$ satisfy suitable regularity and decay properties, with
$\widehat{u}_+$ suitably vanishing on $C_r$ and $\parallel \widehat{u}_+ \parallel_{\infty}$ sufficiently small. Then\par

(1) There exists $T$, $1 \leq T < \infty$, and there exists a unique
solution $(u,A)$ of the (KGS)$_2$ system (\ref{5.1e}) such that $(v, B)
\equiv (u - u_a, A-A_a) \in X([T, \infty ))$.\par

(2) There exists $T$, $1 \leq T < \infty$, and there exists a unique
solution $(u,A)$ of the (KGS)$_2$ system (\ref{5.1e}) such that $(u-u_0, A-A_0)
\in X([T, \infty ))$. One can take the same $T$ and one obtains the same solution as in Part (1).}\\

One should obtain similar results with $X(\cdot )$ replaced by
$X_1(\cdot )$ or by $X_2(\cdot )$, defined by (\ref{5.13e})
(\ref{5.14e}), with suitably stronger regularity and decay assumptions
on $(u_+, A_+, \dot{A}_+)$ and suitably stronger vanishing of $\widehat{u}_+$ on $C_r$. The smallness condition should be the same
in all cases, coming from the contribution of $u_1$ to (\ref{5.6e}),
namely smallness of $\parallel \widehat{u}_+ \parallel_{\infty}$. In particular there
should be no smallness condition on $(A_+, \dot{A}_+)$. \par

It remains to determine appropriate assumptions on
$(u_+, A_+, \dot{A}_+)$ to make the previous pseudoproposition into a
genuine proposition and to remove the conditional form in the
accompanying comments. We refer to \cite{29r} for a set of sufficient
assumptions. On the other hand, as mentioned above, since we are
treating Step 1 in larger spaces than is done in \cite{29r}, those
assumptions are stronger than needed here, and we leave it as an
open question to determine the most economical assumptions adapted to
the present situation and to the various choices of $X(\cdot )$.

We finally mention that a stronger decay, namely $t^{-\lambda}$ with $1
< \lambda <2$, has been obtained in \cite{30r} \cite{31r} on a smaller
space of asymptotic states by the use of more precise asymptotic forms
for ($u_a, A_a)$.\\

\noi{\large\bf Acknowledgements} 

One of us (J. G.) is
grateful to Professor Ozawa and Professor Tsutsumi for the invitation
to the COE Symposium on Nonlinear Dispersive Equations and for the
opportunity to lecture there on the subject matter of this paper. He is
also grateful to Dr. Shimomura for correspondence.

\end{document}